\documentclass[a4paper,10pt]{article}
\usepackage[utf8]{inputenc}
\usepackage[english]{babel}
\usepackage{latexsym,cite,mathrsfs}
\usepackage{amsfonts}
\usepackage{amssymb, amsthm}
\usepackage{xcolor}
\usepackage{marginnote}
\usepackage{amsmath}
\usepackage{a4wide}
\usepackage[affil-it]{authblk}

\numberwithin{equation}{section}

\usepackage{color,enumitem,graphicx}
\usepackage[colorlinks=true,urlcolor=blue,
citecolor=red,linkcolor=blue,linktocpage,pdfpagelabels,
bookmarksnumbered,bookmarksopen]{hyperref}
\usepackage{comment}

\usepackage[left=2.9cm,right=2.9cm,top=2.8cm,bottom=2.8cm]{geometry}
\usepackage[hyperpageref]{backref}

\makeatletter
\providecommand\@dotsep{5}
\def\listtodoname{List of Todos}
\def\listoftodos{\@starttoc{tdo}\listtodoname}
\makeatother

\usepackage{verbatim} 
\usepackage{a4wide}
\usepackage[hyperpageref]{backref}

\usepackage[hyperpageref]{backref}

\newtheorem{thm}{Theorem}[section]
\newtheorem{lem}[thm]{Lemma}
\newtheorem{prop}[thm]{Proposition}
\newtheorem{cor}[thm]{Corollary}
\theoremstyle{definition}
\newtheorem{rem}[thm]{Remark}

\begin{document}

\title{Separating of critical points on the Nehari manifold \\ via the nonlinear generalized Rayleigh quotients}

\author[1]{Marcos L. M. Carvalho\thanks{E-mail: \texttt{marcos$\_$leandro$\_$carvalho@ufg.br}}}
\author[2]{Y.Sh.~Ilyasov \thanks{E-mail: \texttt{ilyasov02@gmail.com}}}
\author[3]{Carlos Alberto Santos\thanks{E-mail: \texttt{csantos@unb.br}}\thanks{Carlos Alberto Santos acknowledges the support of CAPES/Brazil Proc.  $N^o$ $2788/2015-02$}}

\affil[1]{\small 
Instituto de Matem\'atica e Estat\'istica.
Universidade Federal de Goi\'as,
74001-970, Goiania, Brazil}

\affil[2]{\small Institute of Mathematics of UFRC RAS,  450008, Ufa, Russia
\newline
Instituto de Matem\'atica e Estat\'istica.
Universidade Federal de Goi\'as,
74001-970, Goiania, Brazil}

\affil[3]{\small  Universidade de Bras\'ilia, Departamento de Matem\'atica, \newline
 70910-900, Bras\'ilia - DF - Brazil}
\date{}


%
%
%
%

\date{$\ $}
\maketitle

\begin{abstract}
In this paper, we deal with equations of variational form which Nahari manifolds can contain more than two different types of critical points.
We introduce a method of separating critical points on the Nahari manifold, based on the use of nonlinear generalized Rayleigh quotients.
The method is illustrated by establishing the existence of positive solutions, ground states and multiplicity results for a two-parameter nonlinear elliptic boundary problem with polynomial nonlinearities.
\end{abstract}



\section{Introduction}\label{Intr}
The paper deals with the existence of solutions for parametrized  problems
\begin{equation}\label{p*}
\tag{$P$}
	\left\{
	\begin{aligned}
		&-\Delta u= f_{\bar{\lambda}}(x,u)~~\mbox{in}~~\Omega, \\
		&~~~u=0~~\mbox{on}~~\partial\Omega,
	\end{aligned}
	\right. 
\end{equation}
where $f_{\bar{\lambda}}:\Omega\times \mathbb{R} \to \mathbb{R}$, $\bar{\lambda} \in \mathbb{R}^m$, $m\geq 1$, $\Omega \subset \mathbb{R}^N$ is a smooth bounded domain. This equation has a variational form with the energy functional given by
$$
\Phi_{\bar{\lambda}}(u)=\frac{1}{2}\int |\nabla u|^2\,dx+\int F_{\bar{\lambda}}(x,u)\,dx,~~u \in W^{1,2}_0:=W^{1,2}_0(\Omega), 
$$
where $F_{\bar{\lambda}}$ stands for the primitive of $f_{\bar{\lambda}}$ with respect to $u$. 

It is well known that the applicability of variational methods quite depends on the geometry of the nonlinearity  of $\Phi_{\bar{\lambda}}$.
The method that often used in practice for analysing of $\Phi_{\bar{\lambda}}$  is an investigation of one-dimensional geometry of the  fibering functions $\phi_{u, \bar{\lambda}}(t):=\Phi_{\bar{\lambda}}(t u)$, $t\in [0,\infty)$, defined for each $u \in W^{1,2}_0$ (see \cite{ Poh}). The fibering approach is naturally related to the Nehari manifold method (NM-method for short) due the Nehari manifold be defined by means of the critical points of $\Phi_{\bar{\lambda}}(t u)$, namely
$$
\mathcal{N}_{\bar{\lambda}}=\{u\in W^{1,2}_0\setminus\{0\}:~\Phi'_{\bar{\lambda}}(u):=\frac{d}{dt}\Phi_{\bar{\lambda}}(t u)|_{t=1}=0\}.
$$
From this point of view, the complexity of the geometry $ \Phi_{\bar{\lambda}}$,  in the general situation (excluding the degenerate cases), can be ranked depending on the number of critical points of the fibered functions $ \Phi_{\bar{\lambda}}(tu) $, $ t \geq 0 $,  $ u \in W^{1,2}_0 $.
Thus, the simplest case is when $\Phi_{\bar{\lambda}}(t u)$ has at most one critical point (with multiplicity taken into account) for any $u \in W^{1,2}_0\setminus 0$ and  $\bar{\lambda} \in \mathbb{R}^m$. In this case, $\mathcal{N}_{\bar{\lambda}}$ is a $C^1$-manifold and any constrained extremal point of $\Phi_{\bar{\lambda}}$ in $\mathcal{N}_{\bar{\lambda}}$ corresponds to a weak solution of \eqref{p*} (see e.g., \cite{ilyaReil}). However, when $\Phi_{\bar{\lambda}}(t u)$  may have two or more critical points on $\mathcal{N}_{\bar{\lambda}_0}$, for some $u \in W^{1,2}_0\setminus 0$ and $\bar{\lambda}_0 \in \mathbb{R}^m$, the problem becomes more complicated, because of the Nehari manifold $\mathcal{N}_{\bar{\lambda}_0}$ contains points $u$ in which $\Phi''_{\bar{\lambda}_0}(u)=0$ may be true, which prevent us to apply the NM-method directly (see e.g., \cite{ilyaReil}).

Thus, in the case of multiplicity of critical points of the fibering functions $\Phi_{\bar{\lambda}}(t u)$, $t\geq 0$,  $u \in W^{1,2}_0\setminus 0$, one arises the problem of finding domains for the parameters $\bar{\lambda} \in \mathbb{R}^m$ in order to provide the sufficient condition to the applicability of the NM-method (see e.g., \cite{ilyaReil}), namely,
\begin{equation}\label{NM}
\tag{$NM$}
\Phi''_{\bar{\lambda}}(u):=\frac{d^2\Phi_{\bar{\lambda}}}{ds^2}(s u)_{|_{s=1}}\neq 0~\mbox{for all }~~ u \in \mathcal{N}_{\bar{\lambda}}.
\end{equation}
The limit points of parameters $\bar{\lambda}$ in which  \eqref{NM} is satisfied are said to be the  \textit{``extreme values of  the NM-method"}. Recently,  a method  to determine the extreme values of  the NM-method   has been introduced in  \cite{ilyaReil}. This method is based on the investigation of the corresponding \textit{nonlinear generalized Rayleigh quotient} (NG-Rayleigh quotient for short) whose critical values correspond to the extreme values of the NM-method. 

The main purpose of this paper is to suggest a general approach to deal with problems like \eqref{p*} whose fibering functions have more than two critical points. To this end, we develop the nonlinear generalized Rayleigh quotient  method to provide us the extreme values of the NM-method in a recursively way.

Let us outline an idea of our approach.  It is understood, the application of Neahri manifold method requires at least $m-1$ parameters to be taken into account, when the corresponding fibering functions $\phi_{u, \bar{\lambda}}(t)$ may have at most $m$ critical values. Notice that this is in accordance with the conception of Arnol'd's hierarchy of degeneracies (see  \cite{Arnold} pp.15-17). Furthermore, by this conception we need to consider surfaces of unremovable degeneracies (the so-called the dividing boundaries for the domains of generic systems  \cite{Arnold}) of codimension of $m-1$ in the space of all systems.
 In general, this means we have to solve a system of $m$ equations obtained from  $\Phi_{\bar{\lambda}}(t u), \Phi'_{\bar{\lambda}}(t u),\ldots \Phi^{(m)}_{\bar{\lambda}}(t u)$ that, in general, is  a rather complicated problem. 

In the present paper, we overcome this difficulty by applying recursively the nonlinear generalized Rayleigh quotient method \cite{ilyaReil}.  A special feature of this method is that it reduces the number of considered critical points to just one. In other words, it reduces the codimension of degeneracies of the problem (see Figures 1, 2 in \cite{ilyaReil}). As a result of this proceeding, we obtain in the last step the simplest variational problem with zero codimension of degeneracies.
\begin{figure}[!ht]
\begin{minipage}[h]{0.49\linewidth}
\center{\includegraphics[scale=0.7]{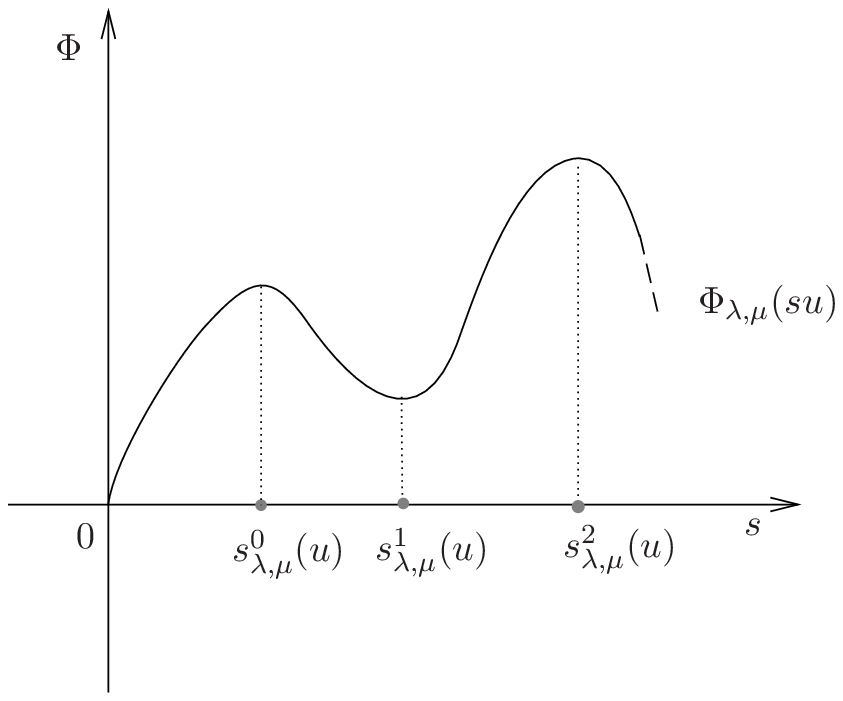}}
\caption{The function $\Phi_{\lambda,\mu}(su)$}
\label{fig1}
\end{minipage}
\hfill
\begin{minipage}[h]{0.49\linewidth}
\center{\includegraphics[scale=0.7]{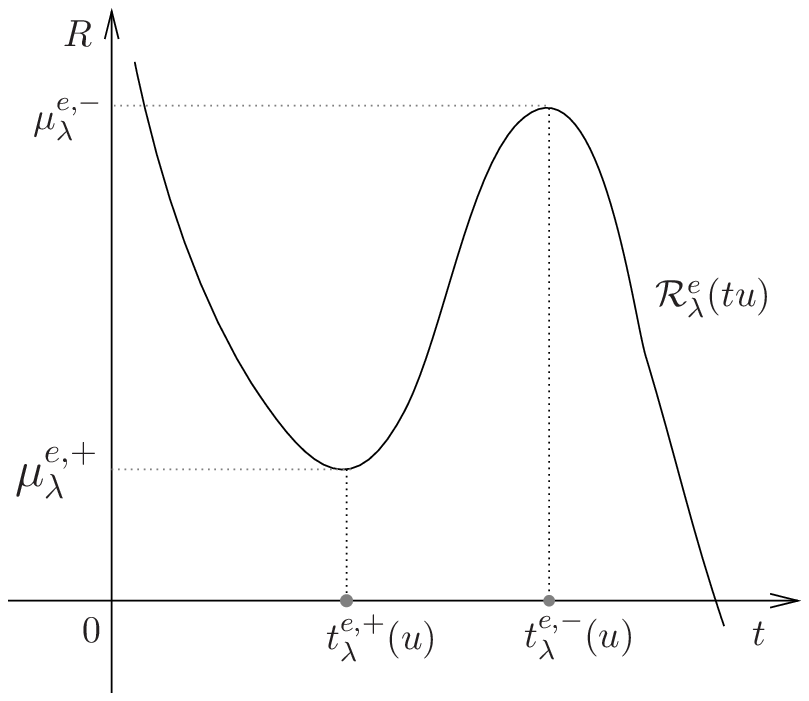}}
\caption{The fibering function $\mathcal{R}^e_\lambda(tu)$}
\label{fig2}
\end{minipage}

\end{figure}

We are going to illustrate the NG-Rayleigh quotient method by applying it to split the Nehari manifold  and find extremal parameters to the 2-parameter ($\bar{\lambda}=(\lambda,\mu)\in \mathbb{R}^2$) model problem
\begin{equation}\label{p}
\left\{
\begin{aligned}
	-&\Delta u= |u|^{\gamma-2}u+\mu|u|^{\alpha-2}u-\lambda |u|^{q-2}u &&\mbox{in}\ \ \Omega, \\
&u=0                                   &&\mbox{on}\ \ \partial\Omega,
\end{aligned}
\right.
\end{equation}
 where $\Omega\subset \mathbb{R}^N$ is a smooth bounded  domain,  $1<q<\alpha<2<\gamma<2^*$ and $\lambda,\mu \in \mathbb{R}$ are parameters. 
Under these assumptions 
the corresponding energy functional
\begin{equation*}
\Phi_{\lambda,\mu}(u)=\frac{1}{2}\int |\nabla u|^2+\frac{\lambda}{q}\int |u|^q-\frac{\mu}{\alpha}\int |u|^\alpha-\frac{1}{\gamma}\int  |u|^\gamma,
\end{equation*}
is well-defined on $W^{1,2}_0$. By a weak solution of \eqref{p} we mean a critical point $u$ of the energy functional $\Phi_{\lambda,\mu}(u)$ on $W^{1,2}_0$.
We are interested in ground states
of \eqref{p}, i.e., a weak solution $u$ of \eqref{p} which satisfies the inequality $\Phi_{\lambda,\mu}(u)\leq \Phi_{\lambda,\mu}(w)$ 
for any non-zero weak solution $w \in 
W^{1,2}_0$ of \eqref{p}.

Observe that the fibering function $\Phi_{\lambda,\mu}(su)$  may have at most three critical points on $s> 0$:
$$0< s^0_{\lambda,\mu}(u)\leq  s^1_{\lambda,\mu}(u) \leq s^2_{\lambda,\mu}(u)<\infty$$ 
which satisfy
$$\Phi''_{\lambda,\mu}(s^0_{\lambda,\mu}(u) u)\leq 0,~~~ \Phi''_{\lambda,\mu}(s^1_{\lambda,\mu}(u) u)\geq 0 ,~~~ \Phi''_{\lambda,\mu}(s^2_{\lambda,\mu}(u) u)\leq 0 .
$$
(see Figure \ref{fig1}). It is not hard to see that, in general,   $s^j_{\lambda,\mu}(u)=s^k_{\lambda,\mu}(u)$ may occur for some $j\neq k$ with $j,k\in \{0,1,2\}$,  and so $\Phi''_{\lambda,\mu}(s^j(u) u)=0$ would be true for some $j\in \{0,1,2\}$, which  may causes difficulties in applying the Nehari manifold method.
An additional difficulty lies in finding a way to separate on the Nehari manifold  the critical points $ s^0 (u) u$ and $ s^2(u)u$, since both of them have the same signs of the derivatives $\Phi''_{\lambda,\mu}$.

To overcome these difficulties, we introduce, in the first step of the recursive procedure,  the \textit{nonlinear generalized Rayleigh quotients}  (see \cite{ilyaReil}) $\mathcal{R}^n_\lambda, \mathcal{R}_\lambda^e: W^{1,2}_0\setminus \{0\} \to \mathbb{R}$   with respect to the parameter $\mu$ 
\begin{eqnarray}\label{RaylQS}
&\mathcal{R}^n_\lambda(u)&=\frac{\int|\nabla u|^2+\lambda\int |u|^q-\int |u|^\gamma}{\int |u|^\alpha}, \\
&\mathcal{R}_\lambda^e(u)&=\frac{\frac{1}{2}\int|\nabla u|^2+\frac{\lambda}{q}\int |u|^q-\frac{1}{\gamma}\int |u|^\gamma}{\frac{1}{\alpha}\int |u|^\alpha},\label{RaylQS2}
\end{eqnarray}
given for $u\in W^{1,2}_0\setminus 0$ and $\lambda\in \mathbb{R}^+$. As consequence, for each $u\in W^{1,2}_0\setminus 0$, we have:
\begin{enumerate}
\item[$(i)$] $\mathcal{R}^n_\lambda(u)=\mu$ if, and only if, $\Phi'_{\lambda,\mu}(u)=0$ and
\item[$(ii)$] $\mathcal{R}^e_\lambda(u)=\mu$ if, and only if, $\Phi_{\lambda,\mu}(u)=0$.
\end{enumerate}
Then following NG-Rayleigh quotient method (see \cite{ilyaReil}), we introduce the so-called \textit{NG-Rayleigh  $\mu$-extremal values} of  $\mathcal{R}^n_\lambda$, $\mathcal{R}^e_\lambda$
given by:
	\begin{align}
&\mu^{e,+}_\lambda=\inf_{u\in  W\setminus 0}\{\mathcal{R}^e_\lambda(u):~(\mathcal{R}_\lambda^e)'(u)=0,~(\mathcal{R}_\lambda^e)''(u)>0\},\label{muplEn}\\
		 &\mu^{e,-}_\lambda=\inf_{u\in  W\setminus 0}\{\mathcal{R}^e_\lambda(u):~(\mathcal{R}_\lambda^e)'(u)=0,~(\mathcal{R}_\lambda^e)''(u)<0\}, \label{muminusEn}\\
		&\mu^{n,+}_\lambda=\inf_{u\in  W\setminus 0}\{\mathcal{R}^n_\lambda(u):~(\mathcal{R}^n_\lambda)'(u)=0,~(\mathcal{R}^n_\lambda)''(u)>0\}, \label{mupl}\\
		&\mu^{n,-}_\lambda=\inf_{u\in  W\setminus 0}\{\mathcal{R}^n_\lambda(u):~(\mathcal{R}^n_\lambda)'(u)=0,~(\mathcal{R}^n_\lambda)''(u)<0\}. \label{muminus}
	\end{align}
In the second step, we  apply the NG-Rayleigh quotient method to the functionals $\mathcal{R}^n_\lambda$ and $\mathcal{R}^e_\lambda$, with respect to the parameter $\lambda$, to obtain the Rayleigh quotients $\Lambda^n,\Lambda^e:W^{1,2}_0\setminus 0 \to \mathbb{R} $ 
\begin{equation}\label{110}
\Lambda^n(u):= \frac{(2-\alpha) \int|\nabla u|^2-(\gamma-\alpha)\int |u|^\gamma}{(\alpha-q)\int |u|^q},
\end{equation}
and 
\begin{equation}\label{111}
\Lambda^e(u):= q\frac{\frac{(2-\alpha)}{2} \int|\nabla u|^2-\frac{(\gamma-\alpha)}{\gamma}\int |u|^\gamma}{(\alpha-q)\int |u|^q},
\end{equation}
which lead to the \textit{NG-Rayleigh $\lambda$-extremal values}  
\begin{align}
&\lambda^{n,*}=\inf_{u\in  W\setminus 0}\sup_{t>0}\Lambda^n(tu),\label{lambExtrem}\\
		 &\lambda^{e,*}=\inf_{u\in  W\setminus 0}\sup_{t>0}\Lambda^e(tu). \label{lambExtrem1}
	\end{align}

Now, we are in position to state our main results. 
The main properties of the extremal values of NG-Rayleigh quotients  are given in the following 
\begin{lem} \label{LemMain}
Assume $1<q<\alpha<2<\gamma<2^*$. Then

	\begin{description}
		\item[(i)] 
		 $0<\lambda^{n,*}<+\infty$, $0<\lambda^{e,*}<+\infty$ and $\lambda^{e,*}<\lambda^{n,*}$.
		\item[(ii)] For any $\lambda \in  (0,\lambda^{e,*})$, there holds  $0\leq\mu^{e,+}_\lambda<\mu^{e,-}_\lambda<+\infty$. Moreover, there exists a minimizer $u^{e,-}_\lambda$ of \eqref{muminusEn}. Furthermore, $u^2_{\lambda, \mu^{e,-}_\lambda}=: u^{e,-}_\lambda$  weakly satisfies equation \eqref{p} with $\mu=\mu^{e,-}_\lambda$ so that  $\Phi_{\lambda,\mu^{e,-}_\lambda}(u^2_{\lambda, \mu^-_\lambda})=0$, $\Phi_{\lambda,\mu^{e,-}_\lambda}''(u^2_{\lambda, \mu^-_\lambda})<0$, $u^2_{\lambda, \mu}\in C^{2}(\Omega)\cap C^{1}(\overline{\Omega})$, $u^2_{\lambda, \mu}>0$ in $\Omega$. 
		\item[(iii)] For any $\lambda \in  (0,\lambda^{n,*})$,  $0\leq \mu^{n,+}_\lambda<\mu^{n,-}_\lambda<+\infty$.
		\item[(iv)] For any $\lambda \in (0,\lambda^{e,*})$, $0\leq  \mu^{n,+}_\lambda\leq \mu^{e,+}_\lambda<\mu^{e,-}_\lambda<\mu^{n,-}_\lambda<+\infty$. 
\end{description}
\end{lem}

%

Our first results on the existence of branch of positive solutions of \eqref{p} is as follows:
\begin{thm}\label{thm2}
Let $1<q<\alpha<2<\gamma<2^*$ and $\lambda \in (0,\lambda^{e,*})$. Assume $\mu \in (\mu^{e,+}_\lambda,\mu^{n,-}_\lambda)$, then problem \eqref{p}
possesses a weak positive solution $u^1_{\lambda, \mu}\in C^{2}(\Omega)\cap C^{1}(\overline{\Omega})$ such that 
$ 
\Phi''_{\lambda,\mu^{e,+}_\lambda}(u^1_{\lambda, \mu})>0,
$
and $\Phi_{\lambda,\mu}(u^1_{\lambda, \mu})<0$.
Moreover, 	$u^1_{\lambda, \mu}$ is a ground state of \eqref{p}.
\end{thm}
\medskip

\noindent The existence of the second branch of positive solutions we prove under the addition restriction $1+\alpha<\gamma<2^*$.

\medskip
\begin{thm}\label{thm3}
	Let $1<q<\alpha<2$, $1+\alpha<\gamma<2^*$ and $\lambda \in (0,\lambda^{n,*})$. Then for any $\mu \in (-\infty, \mu^{n,-}_\lambda)$, problem \eqref{p}
	possesses a weak positive solution $u^2_{\lambda, \mu}\in C^{2}(\Omega)\cap C^{1}(\overline{\Omega})$ such that 
	$$
	\Phi''_{\lambda,\mu}(u^2_{\lambda, \mu})<0, ~~(\mathcal{R}^n_\lambda)''(u^2_{\lambda, \mu})<0.
	$$
	Furthermore, 
	\begin{equation}\label{Phi-}
	\begin{aligned}
		&\Phi_{\lambda,\mu}(u^2_{\lambda, \mu})=0, &&\mbox{if} ~~ \mu=\mu^{e,-}_\lambda,\\
	&\Phi_{\lambda,\mu}(u^2_{\lambda, \mu})<0, &&\mbox{if} ~~\mu \in (\mu^{e,-}_\lambda,\mu^{n,-}_\lambda),\\
	&\Phi_{\lambda,\mu}(u^2_{\lambda, \mu})>0, &&\mbox{if} ~~\mu \in (-\infty, \mu^{e,-}_\lambda).
	\end{aligned}
	\end{equation}
	Moreover, if $\lambda \in (0,\lambda^{e,*})$, $\mu \in (-\infty, \mu^{e,+}_\lambda)$, then	$u^2_{\lambda, \mu}$ is a ground state of \eqref{p}.

\end{thm}

From the just above theorems, we have the following result on the existence of multiple positive solutions.
\begin{cor} Assume $1<q<\alpha<2$ and $1+\alpha<\gamma<2^*$ hold. Then  problem \eqref{p}
admits at least  two distinct positive solutions  $u^1_{\lambda, \mu}$ and $u^2_{\lambda, \mu}$ for any $ 0<\lambda<\lambda^{e,*}$ and $\mu^{e,+}_\lambda<\mu<\mu^{n,-}_\lambda$.

 Furthermore, in the case $1<q<\alpha<2<\gamma<2^*$, problem \eqref{p} has also two distinct positive solutions  $u^1_{\lambda, \mu^{e,-}_\lambda}$ and $u^2_{\lambda, \mu^{e,-}_\lambda}$ for any $ 0<\lambda<\lambda^{e,*}$ and $\mu=\mu^{e,-}_\lambda$.
\end{cor}

It is worth noting that along with the above results, their proof also includes novation. Indeed, below in our proofs  the functionals $\mathcal{R}^n_\lambda(u)$, $\mathcal{R}^e_\lambda(u)$ play essentially role, whose geometric properties, as noted above, are simpler than those of $\Phi_{\lambda,\mu}(u)$.

Notice that problem \eqref{p} in the case $\mu=0$, $\lambda>0$ can not has a non-zero solution $u_\lambda$ such that $\Phi''_{\lambda}(u_{\lambda})\geq 0$. Indeed, for any $u \in W^{1,2}_0\setminus 0 $, the corresponding fibering function $\Phi_{\lambda}(tu)$ has only critical value $t_{max}(u)>0$ with $\Phi''_{\lambda}(t_{max}(u) u)< 0$. This means that it is hardly possible to construct a branch of solutions $u^1_{\lambda, \mu}$ obtained in Theorem \ref{thm2} by local investigation in the neighborhood of point $\mu=0, u=0$, for example, as a bifurcation from zero or by using a priori estimates.
Notice that this implies a conjecture $\mu^{e,+}_\lambda>0$. Moreover, in this connection, the question arises whether it is possible to obtain solutions like $u^1_{\lambda, \mu}$ in Theorem \ref{thm2} without finding the extreme values similar to that of $\mu^{e,+}_\lambda$ and $\mu^{n,-}_\lambda$.

Let us highlight some contribution of our results to literature:
\begin{enumerate}
\item[$(i)$] It is natural to deal with problem \eqref{p} by the Mountain Pass Theorem, introduced in the famous work of Ambrosetty, Rabinowitz \cite{AmbrRabin} due to the geometry of the energy functional. Notice that  the nonlinearity $f_{\mu, \lambda}(u):=|u|^{\gamma-2}u+\mu|u|^{\alpha-2}u-\lambda |u|^{q-2}$  in right hand side of equation \eqref{p} satisfies to the Ambrosetty-Rabinowitz condition \cite{AmbrRabin} for any $\mu, \lambda$, i.e., $0< F_{\mu, \lambda}(s):=\int_0^sf_{\mu, \lambda}(s)\,ds <\theta f_{\mu, \lambda}(s)s$, $\forall |s|>s_0$ for some $s_0>0$ and $\theta \in (1/\gamma,1/2)$.  However, it is not clear that the qualitative information on the solutions derived in above theorems be possible.
\item[$(ii)$]  By this method, it is possible obtaining multiplicity of solutions by separating the Nehari manifolds in different regions given by properties of the  NG-Rayleigh quotient.
\item[$(iv)$] The knowledge on the signs of $\Phi''_{\lambda,\mu}(u_{\lambda, \mu})$, $(\mathcal{R}^n_\lambda)''(u_{\lambda, \mu})$ may permit, in the next steps, to investigate the (in-)stability  of the obtained above solutions for the corresponding non-stationary problem (see e.g., \cite{DiazHYa, payne}) and to get their specific and quality properties  (see e.g., \cite{BobDrIl, DiazHerYA}).
\end{enumerate}

\begin{rem}
It follows from our results  that for any $\lambda \in (0,\lambda^{n,*})$ and $\mu\in (-\infty,\mu^{n,+}_\lambda)$, problem \eqref{p} has no solution $u_{\lambda,\mu}$ such that 
$\Phi''_{\lambda,\mu}(u_{\lambda, \mu})>0$ or  $\Phi''_{\lambda,\mu}(u_{\lambda, \mu})\leq 0$, $\mathcal{R}''_\lambda(u_{\lambda, \mu})\geq 0$ hold.
\end{rem}
\begin{rem}
In the cases, $\lambda=0$, $\mu>0$ or $\mu=0$, $\lambda<0$  problem \eqref{p} coincides with the so-called convex-concave problem first studied in \cite{AmbBrezCer}. The corresponding extremal  values of Nehari manifold method to this type of problems has been studied in \cite{ilConC}.
\end{rem}

\begin{rem}
	In the present paper, we apply the recursive procedure using NG-Rayleigh quotient  method first by the parameter $\mu$ and then by the parameter $\lambda$.  However, this procedure can be changed to the opposite, namely, first by $\lambda$ and then by $\mu$. It seems that the results will be similar to the presented above.  However, this requires additional investigation which we intend to carry out in forthcoming works.
\end{rem}

The paper is organized in the following way. In Section 2, we introduce the  NG-Rayleigh quotients and derive their properties. In Section 3, we present the proofs of Lemma \ref{LemMain}. In Sections 4 and 5, we prove Theorems  \ref{thm2} and \ref{thm3}, respectively. In the Appendix, we present some abstract results from the NG-Rayleigh quotients theory.

\section{ NG-Rayleigh quotients}\label{Sec:2}
In this section, we introduce the   NG-Rayleigh quotients and derive their properties.  We are going to work on the Sobolev space $W^{1,2}_0:=W^{1,2}_0(\Omega)$ endowed with the norm
\begin{equation*}
||u||_{1}=\left( \int_{\Omega }|\nabla u|^{2}\,dx\right) ^{1/2}
\end{equation*}%
and on the Lebesgue spaces $L^{p}(\Omega)$, with $1<p<\infty$, endowed with the standard norm $\|u\|_{L^{p}}$. 

In the first step, we consider the Rayleigh quotients $\mathcal{R}^e_\lambda$,  $\mathcal{R}^n_\lambda$ (see \eqref{RaylQS}-\eqref{RaylQS2}) and  the corresponding fibering Rayleigh quotients
\begin{equation}\label{d21}
\mathcal{R}^e_\lambda(tu):=\frac{\alpha}{\|u\|_{L^{\alpha}}^\alpha}\left(\frac{t^{2-\alpha}}{2}\| u\|^2_1+\frac{t^{q-\alpha}}{q}\lambda\|u\|^q_{L^q}-\frac{t^{\gamma-\alpha}}{\gamma}\|u\|^\gamma_{L^\gamma}\right),~t\geq 0,
\end{equation}
\begin{equation}\label{d22}
\mathcal{R}^n_\lambda(tu):=\frac{1}{\|u\|_{L^{\alpha}}^\alpha}\left(t^{2-\alpha}\| u\|^2_1+t^{q-\alpha}\lambda\|u\|^q_{L^q}-t^{\gamma-\alpha}\|u\|^\gamma_{L^\gamma}\right),~t\geq 0,
\end{equation}
defined for every $u \in W^{1,2}_0\setminus 0$.

Let $\lambda\in \mathbb{R}^+$ and $u \in W^{1,2}_0\setminus 0$.
Simple analysis  shows that the map $\mathcal{R}^e_\lambda(tu)$ ($ \mathcal{R}^n_\lambda(tu)$) may has at most two critical points $0<t_\lambda^{e,+}(u)\leq t_\lambda^{e,-}(u)<+\infty$ ($0<t_\lambda^{n,+}(u)\leq t_\lambda^{n,-}(u)<+\infty$) so that $t_\lambda^{e,+}(u)$  ($t_\lambda^{n,+}(u)$) is a local minimum,  whereas $t_\lambda^{e,-}(u)$ ($t_\lambda^{e,-}(u)$) is a local maximum point of $\mathcal{R}^e_\lambda(tu)$ ($ \mathcal{R}^n_\lambda(tu)$) (see Figure \ref{fig2}). 
We are interested in finding the values $\lambda>0$ for which the function $\mathcal{R}^e_\lambda(tu)$ ($\mathcal{R}^n_\lambda(tu)$) has precisely two distinct critical points $0<t_\lambda^{e,+}(u)<t_\lambda^{e,-}(u)<+\infty$ ($0<t_\lambda^{n,+}(u)<t_\lambda^{n,-}(u)<+\infty$) such that $(\mathcal{R}^e_\lambda)''(t_\lambda^{e,\pm}(u)u)\neq 0$  ($(\mathcal{R}^n_\lambda)''(t_\lambda^{n,\pm}(u)u)\neq 0$) for any $u \in W^{1,2}_0\setminus 0$. To this end, following NG-Rayleigh quotient method \cite{ilyaReil}, we  consider the equation
\begin{align*}
	&(\mathcal{R}^e_\lambda)'(tu)\equiv ~~\frac{\alpha}{t\|u\|_{L^{\alpha}}^\alpha} \times \\
	&\left(\frac{(2-\alpha)t^{2-\alpha}}{2}\| u\|^2_1+\frac{(q-\alpha)t^{q-\alpha}}{q}\lambda\|u\|^q_{L^q}-\frac{(\gamma-\alpha)t^{\gamma-\alpha}}{\gamma}\|u\|^\gamma_{L^\gamma}\right)=0.
\end{align*}
Solving this  with respect to $\lambda$ we obtain the fibering Rayleigh quotient (cf. \eqref{lambExtrem})
\begin{equation}\label{LambdaE}
\Lambda^e(tu):= \frac{q}{(\alpha-q)\|u\|_{L^q}^q}\left(t^{2-q}\frac{(2-\alpha)}{2} \| u\|^2_1-t^{\gamma-q}\frac{(\gamma-\alpha)}{\gamma}\|u\|^\gamma_{L^\gamma}\right),
\end{equation}
given for $t>0$ and $u \in W^{1,2}_0\setminus 0$. In a similar way, we obtain 
\begin{equation}
\Lambda^n(tu):= \frac{1}{(\alpha-q)\|u\|_{L^q}^q}\left((2-\alpha) t^{2-q}\| u\|^2_1-(\gamma-\alpha)t^{\gamma-q}\|u\|^\gamma_{L^\gamma}\right),
\end{equation}
for $t>0,$ and $u \in W^{1,2}_0\setminus 0$.
Let us note that finding critical points of the fibering functions of $\mathcal{R}^e_\lambda(tu)$ and $\mathcal{R}^n_\lambda(tu)$ for  $u \in W_0^{1,2}$, are equivalent to solving the equations $\Lambda^e(tu)=\lambda$ and $\Lambda^n(tu)=\lambda$, respectively. 

\begin{figure}[!ht]
\begin{minipage}[h]{0.49\linewidth}
\center{\includegraphics[scale=0.7]{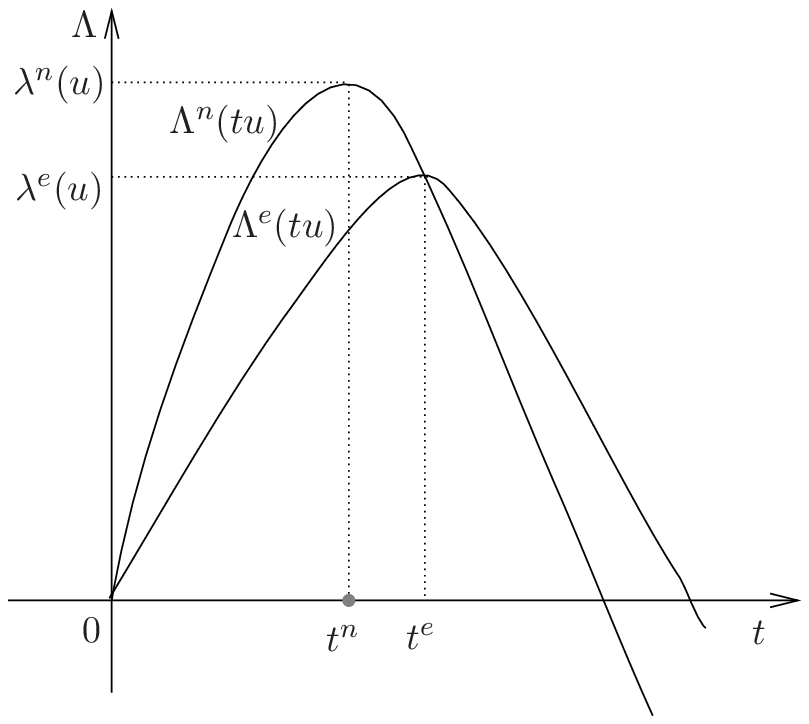}}
\caption{The functions $\mathcal{R}^e_\lambda(tu)$, $\mathcal{R}^n_\lambda(tu)$}
\label{fig3}
\end{minipage}
\hfill
\begin{minipage}[h]{0.49\linewidth}
\center{\includegraphics[scale=0.7]{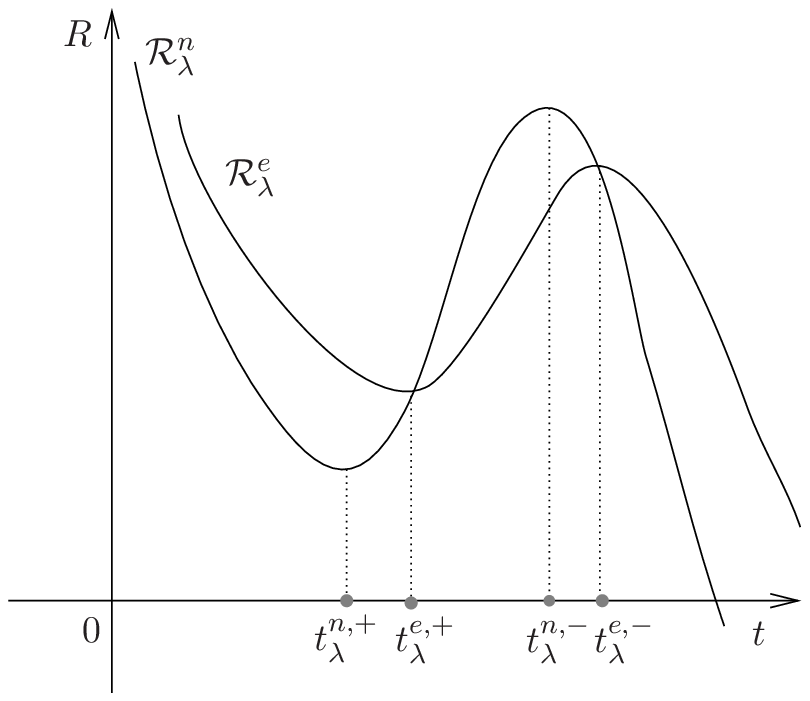}}
\caption{The functions $\Lambda^n(tu)$, $\Lambda^e(tu)$}
\label{fig4}
\end{minipage}
\end{figure}

It is easily to see that  for any $u\in  W^{1,2}_0\setminus 0$, the function $\Lambda^e(tu)$ has a unique critical point $t^e(u)>0$ which is a global maximum point of $\Lambda^e(tu)$ (see Figure \ref{fig3}). Moreover,  
$\Lambda^e(tu)|_{t=0}=0$ and $\Lambda^e(tu) \to -\infty$ as $t \to +\infty$.
Solving the equation 
$
\frac{d}{dt}\Lambda^e(tu)=0
$
we find 
$$
t^e(u)= \left(C_e\frac{\| u\|^2_1}{\| u\|^\gamma_{L^\gamma}}\right)^{1/(\gamma-2)},
$$ 
where $ C_e=\frac{\gamma(2-\alpha)(2-q)}{2(\gamma-\alpha)(\gamma-q)}$. Thus we have the following  \textit{nonlinear generalized Rayleigh $\lambda$-quotient} (the\textit{ NG-Rayleigh $\lambda$-quotient} for short) (see \cite{ilyaReil}) 
$$
\lambda^{e,*}(u):=\Lambda^e(t^e(u)u)=c^e_{q,\gamma}\frac{\| u\|_1^{2\frac{\gamma-q}{\gamma-2}}}{\| u\|^q_{L^q} \cdot\| u\|_{L^\gamma}^{\gamma\frac{2-q}{\gamma-2}}},
$$
where 
\begin{equation}\label{const}
	c^e_{q,\gamma}=\frac{q\gamma^{\frac{2-q}{\gamma-2}}}{2^{\frac{\gamma-q}{\gamma-2}}}\cdot \frac{(2-\alpha)^{\frac{\gamma-q}{\gamma-2}}(2-q)^{\frac{2-q}{\gamma-q}}(\gamma-2)}{(\alpha-q)(\gamma-\alpha)^{\frac{2-q}{\gamma-2}}(\gamma-q)^{\frac{\gamma-q}{\gamma-2}}}.
\end{equation}
Similarly the function $\Lambda^n(tu)$ achieves a global maximum and  the corresponding maximum point $t^n(u)$ can be found in the explicit form $t^n(u)=(C_n(\alpha,q,\gamma)\| u\|^2_1/\| u\|^\gamma_{L^\gamma})^{1/(\gamma-2)}$, where $C_n(\alpha,q,\gamma)=\frac{(2-\alpha)(2-q)}{(\gamma-\alpha)(\gamma-q)}$,
Hence we have the following \textit{NG-Rayleigh $\lambda$-quotient }
$$
\lambda^n(u):=\Lambda^n(t^n(u)u)= c^n_{q,\gamma}\frac{\| u\|_1^{2\frac{\gamma-q}{\gamma-2}}}{\| u\|^q_{L^q} \cdot\| u\|_{L^\gamma}^{\gamma\frac{2-q}{\gamma-2}}},
$$
where 
\begin{equation}\label{const2}
	c^n_{q,\gamma}=\frac{2^{\frac{\gamma-q}{\gamma-2}}}{q\gamma^{\frac{2-q}{\gamma-2}}}c^e_{q,\gamma}.
\end{equation}
As a result, we have the \textit{NG-Rayleigh  $\lambda$-extremal values} (cf. \eqref{lambExtrem})
\begin{align}
	&\lambda^{e,*}=\inf_{u\in  W^{1,2}_0\setminus 0}\lambda^{e,*}(u)\equiv \inf_{u\in  W^{1,2}_0\setminus 0} \sup_{t>0} \Lambda^e(tu),\\
		&\lambda^{n,*}=\inf_{u\in  W^{1,2}_0\setminus 0}\lambda^{n,*}(u)\equiv \inf_{u\in  W^{1,2}_0\setminus 0} \sup_{t>0} \Lambda^n(tu).
\end{align}
It is not hard to show using Holder's and Soblev's inequalities that
\begin{align} \label{eq:lambdaEB}
	0<\lambda^{e,*}<+\infty,
	~~0<\lambda^{n,*}<+\infty.
\end{align}

Let us now introduce the  NG-Rayleigh $\mu$-quotients corresponding to $\mathcal{R}^n_\lambda$ and $\mathcal{R}^e_\lambda$. To this end, we need
\begin{prop}\label{propRE} 
	For each $\lambda\in (0,\lambda^{e,*})$ and $u\in W^{1,2}_0\setminus 0$, the function $\mathcal{R}^e_\lambda(tu)$ has  two distinct critical points such that $0<t_\lambda^{e,+}(u)<t_\lambda^{e,-}(u)$. Moreover, 
		\par (i)\, $t_\lambda^{e,+}(u)$ is a local minimum point such that $(\mathcal{R}^e_\lambda)''(t_\lambda^{e,+}(u)u)>0$ holds strongly and $t_\lambda^{e,-}(u)$ is a local maximum such that $(\mathcal{R}^e_\lambda)''(t_\lambda^{e,-}(u)u)<0$ holds strongly;
		\par (ii)\, $t_\lambda^{e,+}(u)$, $t_\lambda^{e,-}(u)$ are $C^1$-functional on $W^{1,2}_0\setminus 0$.

\end{prop}

	\begin{proof} 	Let $\lambda\in (0,\lambda^{e,*})$ and $u\in W^{1,2}_0\setminus 0$. 
	Then in view of \eqref{lambExtrem1}, $\mathcal{R}(t_e(u)u)>\lambda$ because $t^e(u)$ is the global maximum point of the function $t\mapsto \Lambda(tu)$. This implies that the equation  $\Lambda^e(tu)=\lambda$ has precisely two roots $t^{e,+}_\lambda(u), t^{e,-}_\lambda(u)$ such that $0<t^{e,+}_\lambda(u)<t^e(u)<t^{e,-}_\lambda(u)$
	and 
		$$
			(\Lambda^e)'(t_\lambda^{e,+}(u)u)>0, ~~ (\Lambda^e)'(t_\lambda^{e,-}(u)u)<0.
	$$
Now applying Proposition \ref{Prop1}  with $\mathcal{R}(tu)$, $\Phi_\nu(tu)$  replaced by $\Lambda^e(tu)$,  $\mathcal{R}^e_\lambda(tu)$, respectively,  we obtain assertion \textit{(i)}.

Since  $\Lambda^e(tu)
\in C^1(\mathbb{R}^+\times W^{1,2}_0)$ and $\frac{d}{dt}\Lambda^e(tu)|_{t=t_\lambda^{e,\pm}(u)}=(\Lambda^e)'(t_\lambda^{e,\pm}(u)u)\neq 0$ for  $\lambda\in (0,\lambda^{e,*})$, the proof of \textit{(ii)} follows by the Implicit Function Theorem. 
	\end{proof}
In a similar way, it can be proven 
\begin{prop}\label{propRN2} 
	For each $\lambda\in (0,\lambda^{n,*})$ and $u\in W^{1,2}_0\setminus 0$, the function $\mathcal{R}^n_\lambda(tu)$ has precisely two distinct critical points such that $0<t_\lambda^{n,+}(u)<t_\lambda^{n,-}(u)$. Moreover, 
		\par (i)\, $t_\lambda^{n,+}(u)$ is a local minimum point such that $(\mathcal{R}^n_\lambda)''(t_\lambda^{n,+}(u)u)>0$ and $t_\lambda^{n,-}(u)$ is a local maximum such that $(\mathcal{R}^n_\lambda)''(t_\lambda^{n,-}(u)u)<0$;
		\par (ii)\, $t_\lambda^{n,+}(u)$, $t_\lambda^{n,-}(u)$ are $C^1$-functional on $W^{1,2}_0\setminus 0$. 
\end{prop}

As a consequence of 
Propositions \ref{propRE}, \ref{propRN2}, we are able to introduce the following \textit{NG-Rayleigh $\mu$-quotients} 
\begin{align}
	\mu^{e,+}_\lambda(u):&=\mathcal{R}^e_\lambda(t_\lambda^{e,+}(u)u), ~~~ \mu^{e,-}_\lambda(u):=\mathcal{R}^e_\lambda(t_\lambda^{e,-}(u)u), \label{mumu}\\
	\mu^{n,+}_\lambda(u):&=\mathcal{R}^n_\lambda(t_\lambda^{n,+}(u)u), ~~~\mu^{n,-}_\lambda(u):=\mathcal{R}^n_\lambda(t_\lambda^{n,-}(u)u), ~~ \label{muMinus}
\end{align}
for $u\in  W^{1,2}_0\setminus 0$. Furthermore, Propositions \ref{propRE}, \ref{propRN2} imply that $0<\mu^{e,+}_\lambda(u)<\mu^{e,-}_\lambda(u)$, ~ $0<\mu^{n,+}_\lambda(u)<\mu^{n,-}_\lambda(u)$ for every $u\in  W^{1,2}_0\setminus\{ 0\}$ and $\mu^{e,+}_\lambda(u),\mu^{e,-}_\lambda(u)$,  $\mu^{n,+}_\lambda(u),\mu^{n,-}_\lambda(u)$  are $C^1$ and $0$-homogeneous functionals on $W^{1,2}_0\setminus 0$.

We need in the following properties of $\mu^{e,+}_\lambda(u)$.
\begin{lem}\label{lemRE}
 Let $\lambda\in (0,\lambda^{e,*})$. Then 
 any critical point $\bar{u} \in  W^{1,2}_0\setminus 0$ of  $\mu^{e,-}_\lambda(u)$  weakly satisfies   \eqref{p} with $\mu=\mu^{e,-}_\lambda(\bar{u})$  and it has a zero energy level, i.e., $\Phi_{\lambda,\mu^{e,-}_\lambda(\bar{u})}(\bar{u})=0$.
\end{lem}
\begin{proof}
 Let $\bar{u} \in  W^{1,2}_0\setminus 0$ be a critical point of $\mu^{e,-}_\lambda(u)$. Then by \eqref{RaylQS2} we have $\Phi_{\lambda,\mu^{e,-}_\lambda(\bar{u})}(\bar{u})=0$, i.e., $\bar{u}$ lies on zero energy level. Now, using (ii), Proposition \ref{propRE}, we derive
\begin{align*}
	0=&D_u\mu^{e,-}_\lambda(\bar{u})(\psi)=D_u\mathcal{R}^e_\lambda(t_\lambda^{e,-}(\bar{u})\bar{u})(\psi)=\\
	&\frac{\partial}{\partial t}\mathcal{R}^e_\lambda(t\bar{u})|_{t=t_\lambda^{e,-}(\bar{u})}(\bar{u}D_ut_\lambda^{e,-}(\bar{u})(\psi))+D_u\mathcal{R}^e_\lambda(t\bar{u})|_{t=t_\lambda^{e,-}(\bar{u})}(t_\lambda^{e,-}(\bar{u})\psi)=\\
	&\frac{D_ut_\lambda^{e,-}(\bar{u})(\psi)}{t_\lambda^{e,-}(\bar{u})}(\mathcal{R}^e_\lambda)'(t_\lambda^{e,-}(\bar{u})\bar{u})+\frac{\alpha t_\lambda^{e,-}(\bar{u})}{\|(t_\lambda^{e,-}(\bar{u})\bar{u}\|_{L^\alpha}^\alpha} D_u\Phi_{\lambda,\mu}(t_\lambda^{e,-}(\bar{u})\bar{u})(\psi),
\end{align*}
$\forall \psi \in W^{1,2}_0$. Since $(\mathcal{R}^e_\lambda)'(t_\lambda^{e,-}(\bar{u})\bar{u})=0$, we obtain that $\bar{u}$ is a critical point of $\Phi_{\lambda,\mu}$, i.e., $\bar{u}$ is a weak solution of \eqref{p}. 
\end{proof}

We need also accounts for the locations of functions $\Lambda^e(tu),\Lambda^n(tu)$ and $\mathcal{R}_\lambda^e(tu)$, $\mathcal{R}_\lambda^n(tu)$ relative to each other.	
	\begin{lem}\label{lemGEOM}
		Let $u\in W^{1,2}_0\setminus 0$. 
		\begin{description}
			\item[(i)] $\Lambda^e(tu)<\Lambda^n(tu)$, for sufficiently small $t>0$;
			\item[(ii)] $\Lambda^{e}(tu)=\Lambda^n(tu),~t>0,$ if and only if, $t=t^e(u)$;
			\item[(iii)] $t_\lambda^{n,+}(u)<t_\lambda^{e,+}(u)<t^e(u)<t_\lambda^{n,-}(u)<t_\lambda^{e,-}(u)$ for each $\lambda \in(0,\lambda^{e,*})$;
			\item[(iv)] for each $\lambda \in(0,\lambda^{e,*})$, the equation $\mathcal{R}_\lambda^e(tu)=\mathcal{R}_\lambda^n(tu)$ has precisely two solutions $t=t_\lambda^{e,+}(u)$ and $t=t_\lambda^{e,-}(u)$.
		\end{description}
	\end{lem}
	
	\begin{proof}
			
				Let $u\in W^{1,2}_0\setminus0$. 
			 Observe, 
				$$
				\lim_{t\rightarrow 0}\frac{\Lambda^e(tu)}{\Lambda^n(tu)}=\lim_{t\rightarrow 0}\frac{q\left(\frac{(2-\alpha)}{2} \| u\|^2_1-t^{\gamma-2}\frac{(\gamma-\alpha)}{\gamma}\|u\|^\gamma_{L^\gamma}\right)}{(2-\alpha) \| u\|^2_1-(\gamma-\alpha)t^{\gamma-2}\|u\|^\gamma_{L^\gamma}}=\frac{q}{2}<1.$$		
Thus we get \textbf{(i) }. The equality
				$\Lambda_{e}(tu)=\Lambda^n(tu)$ is equivalent to
				$$(2-\alpha) t^{2-q}\| u\|^2_1-(\gamma-\alpha)t^{\gamma-q}\|u\|^\gamma_{L^\gamma}=q\left(t^{2-q}\frac{(2-\alpha)}{2} \| u\|^2_1-t^{\gamma-q}\frac{(\gamma-\alpha)}{\gamma}\|u\|^\gamma_{L^\gamma}\right).$$
This implies
				$$
				0=\frac{(2-q)(2-\alpha)}{2} t^{1-q}\| u\|^2_1-\frac{(\gamma-q)(\gamma-\alpha)}{\gamma}t^{\gamma-q-1}\|u\|^\gamma_{L^\gamma}=(\Lambda^{e}(tu))',
				$$
which gives \textbf{(ii)}. Proof of \textbf{(iii)} directly follows from  \textbf{(i),(ii)}.

			Observe, 	$\mathcal{R}_\lambda ^e(tu)=\mathcal{R}_\lambda^n(tu)$ implies 
				$$t^{2-\alpha}\| u\|^2_1+\lambda t^{q-\alpha}\|u\|^q_{L^q}-t^{\gamma-\alpha}\|u\|^\gamma_{L^\gamma}=\frac{\alpha t^{2-\alpha}}{2}\| u\|^2_1+\frac{\lambda\alpha t^{q-\alpha}}{q}-\frac{\alpha t^{\gamma-\alpha}}{\gamma}\|u\|^\gamma_{L^\gamma}.$$
				Hence,
						\begin{align*}
					0=\frac{(2-\alpha)}{2}t^{2-\alpha}\| u\|^2_1-\frac{\gamma-\alpha}{\gamma}t^{\gamma-\alpha}\|u\|^\gamma_{L^\gamma}-&\frac{\lambda(\alpha-q)}{q}t^{q-\alpha}\|u\|^q_{L^q}=\\
					&\frac{(\alpha-q)\|u\|^q_{L^q}t^{q-\alpha}}{q}\left(\Lambda_{e}(tu)-\lambda\right).
				\end{align*}
		Thus we get		\textbf{(iv)}.

\end{proof}
Summarizing the above information, we have the situation described in the Figures \ref{fig3},\ref{fig4}.
	
Below, we give more refined information about the convexity of the NG-Rayleigh quotient of $\mathcal{R}^n_\lambda$.

\begin{prop}\label{prop:2deriv}
 Assume 	$1<q<\alpha<2$, and  $1+\alpha<\gamma<2^*$ hold. Then there exists a one root $r^n_\lambda(u) \in (t_\lambda^{n,+}(u), t_\lambda^{n,-}(u))$ so that  $(\mathcal{R}^n_\lambda)''(tu)>0$ for $t \in (0,r^n_\lambda(u))$, $(\mathcal{R}^n_\lambda)''(r^n_\lambda(u)u)=0$ and $(\mathcal{R}^n_\lambda)''(tu)<0$ for $t \in (r^n_\lambda(u), +\infty)$ for any  $u \in W_0^{1,2}\setminus\{0\}$ and $0 < \lambda <\lambda^{n,*}$
\end{prop}
\begin{proof}
The claim follows from the expression 
	\begin{align*}
		(\mathcal{R}_{\lambda}^n)''(tu)&=\frac{t^{-\alpha}}{\|u\|_{L^{\alpha}}^\alpha}\big[(2-\alpha)(1-\alpha)\|u\|^2_1+(q-\alpha)(q-\alpha-1)\lambda \|u\|_{L^q}^qt^{q-2}\\
		&-(\gamma-\alpha)(\gamma-\alpha-1)\lambda \|u\|_{L^\gamma}^\gamma t^{\gamma-2}\big]
	\end{align*}
and the assumption $2<1+\alpha<\gamma$.
\end{proof}

To end this section, let us turn our attention to solutions of the equations
$$
\mathcal{R}^n_\lambda(tu)= \mu~~\mbox{and}~~ \mathcal{R}^e_\lambda(tu)=\mu, ~~ u \in W_0^{1,2}\setminus\{0\},\mu \in \mathbb{R}
$$
We know from above information that the number of solutions to these equations depend on the value of the parameters $\lambda>0$ and $\mu\in \mathbb{R}$. Let us understand this claim to $\mathcal{R}^n_\lambda(tu)$, $t\geq 0$. 

Let $u \in W_0^{1,2}\setminus\{0\}$. So, it follows from Lemma \ref{lemGEOM} and \eqref{lambExtrem}, \eqref{muMinus} that:
\begin{enumerate}
\item[$(i)$] $\Lambda^n(t u) =\lambda$ has two different solutions if $\lambda< \lambda^{n,*}(u)$, that is, the fibering function $\mathcal{R}_{\lambda}^n(tu)$, $t> 0$, has two critical points. This implies that equation $\mathcal{R}_{\lambda}^n(tu)=\mu$ has:
\begin{enumerate}
\item[$(i)_1$] three different solutions, say: 
$$0<s^0_{\lambda,\mu}(u)< s^1_{\lambda,\mu}(u) <s^2_{\lambda,\mu}(u)<\infty$$
if $\mu^{n,+}_\lambda(u) <\mu< \mu^{n,-}_\lambda(u)$;
\item[$(i)_2$] two solutions if $\mu = \mu^{n,+}_\lambda(u)$, in this case, $s^0_{\lambda,\mu}(u) = s^1_{\lambda,\mu}(u) = t_\lambda^{n,+}(u)$; 
\item[$(i)_3$] two solutions if $\mu = \mu^{n,-}_\lambda(u)$, in this case, $s^1_{\lambda,\mu}(u) = s^2_{\lambda,\mu}(u) = t_\lambda^{n,-}(u)$; 
\item[$(i)_4   $] one solution if either   $\mu <\mu^{n,+}_\lambda(u) $ (just $s^0_{\lambda,\mu}(u)$) or  $\mu >  \mu^{n,-}_\lambda(u)$ $s^2_{\lambda,\mu}(u)$);
\end{enumerate}
\item[$(ii)$] $\Lambda^n(t u) =\lambda$ has just one solution if $\lambda= \lambda^{n,*}(u)$, that is, $\mathcal{R}_{\lambda}^n(tu)$, $t\geq 0$, has an only one critical point $ t_\lambda^{n,+}(u)= t_\lambda^{n,-}(u)$ and  $\mathcal{R}_{\lambda}^n(tu)=\mu$ has an only  solution for any $\mu \in \mathbb{R}$,
\item[$(iii)$] $\Lambda^n(t u) =\lambda$ has no solution if $\lambda> \lambda^{n,*}(u)$, that is, $\mathcal{R}_{\lambda}^n(tu)$, $t\geq 0$, has no critical point and  $\mathcal{R}_{\lambda}^n(tu)=\mu$ has an only  solution for any $\mu \in \mathbb{R}$ as well.
\end{enumerate}

All in all, we have.
\begin{prop}\label{p28} Assume $0<\lambda<\lambda^{n,*}$ and $u\in W^{1,2}_0\setminus 0$. If $\mu^{n,+}_\lambda(u) <\mu< \mu^{n,-}_\lambda(u)$, then the equation $\mathcal{R}^n_\lambda(tu)=\mu$, $t> 0$ has:
\begin{enumerate}
\item[$(i)$] three solutions $s^0_{\lambda,\mu}(u),s^1_{\lambda,\mu}(u),s^2_{\lambda,\mu}(u)$ such that 
$$0<s^0_{\lambda,\mu}(u)< t_\lambda^{n,+}(u)<s^1_{\lambda,\mu}(u) <t_\lambda^{n,-}(u)<s^2_{\lambda,\mu}(u)<\infty$$
and
$$(\mathcal{R}_{\lambda}^n)'(s^0_{\lambda,\mu}(u)u)<0,~~(\mathcal{R}_{\lambda}^n)'(s^1_{\lambda,\mu}(u)u)>0~\mbox{and}~ (\mathcal{R}_{\lambda}^n)'(s^2_{\lambda,\mu}(u)u)<0,$$
\item[$(ii)$]  $s^j_{\lambda,\mu}:W_0^{1,2}\setminus\{0\} \to (0,\infty)$ is $C^1$- functional, for $j=0,1,2$;
\end{enumerate}
\end{prop}
\begin{proof} The proof of item $(i)$ follows from the above information. The proof of item $(ii)$ follows from arguments like those done to prove $(ii)$, Proposition \ref{propRE}. This ends the proof.
\end{proof}

\section{ Proof of Lemma \ref{LemMain}}\label{Sec:3}
By Propositions \ref{propRE}, \ref{propRN2}, the NG-Rayleigh extremal values \eqref{muplEn}-\eqref{muminus} are equivalent to the definitions \eqref{mumu} and \eqref{muMinus}, respectively. Furthermore, 
evidently that $0\leq \mu^{e,+}_\lambda,\mu^{e,-}_\lambda,\mu^{n,+}_\lambda,\mu^{n,-}_\lambda<+\infty$ for $\lambda \in (0,\lambda^{e,*})$ and $\lambda \in (0,\lambda^{n,*})$, respectively. 

\begin{prop}\label{lem22}
\begin{description}
	\item[(i)] For any $\lambda \in (0,\lambda^{e,*})$, variational problem  \eqref{muminusEn} has a minimizer $ u^{e,-}_\lambda \in W^{1,2}_0\setminus 0$ so that 
 $0<\mu^{e,-}_\lambda=\mu^{e,-}_\lambda(u^{e,-}_\lambda)$.
\item[(ii)] For any $\lambda \in (0,\lambda^{n,*})$ variational problems  \eqref{muminus} has a minimizer $ u^{n,-}_\lambda \in W^{1,2}_0\setminus 0$  so that 
 $0<\mu^{n,-}_\lambda=\mu^{n,-}_\lambda(u^{n,-}_\lambda)$.
\end{description}
\end{prop}
		\begin{proof} \textbf{(i)} 
Observe that for any sequences  $(u_m)\subset W^{1,2}_0\setminus 0$ the assumption $t_\lambda^{e,-}(u_m)\to \infty$  entails $\mu^{e,-}_\lambda(u_m)\to \infty$. Indeed, since $\mu^{e,-}_\lambda(u)$ is a homogeneous functional we 
	may assume that $\|u_m\|=1$. Hence, since  $(\mathcal{R}^e_\lambda)'(t_\lambda^{e,-}(u_m)u_m) =0$, we have
		\begin{equation}\label{RSCHT}
			\|u_m\|_{L^\gamma}^\gamma=\frac{\gamma(2-\alpha)}{2(\gamma-\alpha)}t_\lambda^{e,-}(u_m)^{2-\gamma}-\lambda\frac{\gamma(\alpha-q)}{q(\gamma-\alpha)}\|u_m\|_{L^q}^qt_\lambda^{e,-}(u_m)^{q-\gamma},
		\end{equation}
and therefore
	
		\begin{align*}
			\lim_{m\to \infty}&\mathcal{R}^e(t_\lambda^{e,-}(u_m)u_m)=\\
			&~~~\lim_{m\to \infty}\frac{\left[\frac{(\gamma-2)}{2}t_\lambda^{e,-}(u_m)^{2-\alpha}+\lambda\frac{(\gamma-q)}{q}\|u_m\|_{L^q}^qt_\lambda^{e,-}(u_m)^{q-\alpha}\right]}{(\gamma-\alpha)\|u_m\|_{L^{\alpha}}^\alpha}=\infty.
		\end{align*}
Denote by $(u_m)\subset W^{1,2}_0$ the minimizer sequence of  \eqref{muminusEn}, i.e.,
	$$
	\mu^{e,-}_\lambda(u_m)=\mathcal{R}^e_\lambda(t_m u_m)\rightarrow \mu^{e,-}_\lambda,
	$$
	where $t_m:=t_\lambda^{e,-}(u_m)$ and we  assume that $\|u_m\|=1$, $m=1,\ldots$. Hence the Sobolev embedding  and  Banach-Alaoglu theorems imply that there exists a subsequence, which we again denote by $(u_m)$, such that $u_m\rightharpoondown u^{e,-}_\lambda$ weakly in $W^{1,2}_0$ and strongly $u_m\to u^{e,-}_\lambda$ in $L^p$, $1\leq p<2^*$.  By the above $(t_m)$ is bounded. Hence, up to a subsequence, we have
	$
	t_m \to \bar{t}~~\mbox{as}~~m \to +\infty
	$,
	for some $\bar{t} \in [0,+\infty)$.
	Notice that by Lemma \ref{lemGEOM} and since $\|u_m\|=1$, we have 
	,  we have
$$
t_m=t_\lambda^{e,-}(u_m)>t^{e}(u)=C^e(\alpha,q,\gamma)\left(\frac{\|u_m\|_1^2}{\|u_m\|_{L^\gamma}^{\gamma}}\right)^{1/(\gamma-2)}\geq C^e(\alpha,q>0,\gamma).
$$
Hence	$\bar{t} \neq 0$. Suppose that $u^{e,-}_\lambda=0$. Then $\|u_m\|_{L^\gamma}, \|u_m\|_{L^q} \to 0$. However, since $t_m \to \bar{t}>0$, this implies, in view of \eqref{RSCHT}, a contradiction. Thus, we have proved that $\mu^{e,-}_\lambda >0$ and $u^{e,-}_\lambda \neq 0$.	

By the weak low  semi-continuity  of $\|\cdot\|_1$, 
$\|\bar{u}^\pm\|_1\leq \liminf_{n\rightarrow\infty}\|u_m\|_1$. Obviously if here the equality holds, then $u^{e,-}_\lambda$ is a minimizer of \eqref{muminusEn} and the proof is completed.
Assume the converse $\|u^{e,-}_\lambda\|< \liminf_{m\rightarrow\infty}\|u_m\|$. Then 
	\begin{align}
		&\tilde{\mu}:=\mathcal{R}^e_\lambda(\bar{t}u^{e,-}_\lambda)<	\liminf_{m \rightarrow\infty}\mathcal{R}^e_\lambda(t_\lambda^{e,-}(u_m)u_m)=\mu^{e,-}_\lambda, \label{eq:comm}\\
		&(\mathcal{R}^e_\lambda)'(\bar{t}u^{e,-}_\lambda)<(\mathcal{R}^e_\lambda)'(t_\lambda^{e,-}(u_m)u_m)=0.\label{eq:comm2}
	\end{align}
Furthermore,  $0=(\mathcal{R}^e_\lambda)'(t_\lambda^{e,-}(\bar{u})\bar{u})<\liminf_{n\rightarrow\infty}(\mathcal{R}^e_\lambda)'(t_\lambda^{e,-}(\bar{u})u_m)$ and therefore for sufficiently large $m$ we have  $(\mathcal{R}^e_\lambda)'(t_\lambda^{e,-}(\bar{u})u_m)>0$. 
 Hence  $t_\lambda^{e,-}(\bar{u})<t_\lambda^{e,-}(u_m)$ for sufficiently large $m$ and thus 
\begin{align*}
	\mu^{e,-}_\lambda(\bar{u})= \mathcal{R}^e_\lambda(t_\lambda^{e,-}(\bar{u})\bar{u})<\liminf_{m\rightarrow\infty}&\mathcal{R}^e_\lambda(t_\lambda^{e,-}(\bar{u})u_m)<\\
	&\liminf_{m\rightarrow\infty}\mathcal{R}^e_\lambda(t_\lambda^{e,-}(u_m)u_m)=\mu^{e,-}_\lambda
\end{align*}
which is a contradiction. 	Thus we get \textbf{(i)}.
	
	The proof of \textbf{(ii)} follows by the same method as in \textbf{(i)}.
\end{proof}
From here it follows 
\begin{cor}\label{cor: ineqStr}
	
	\begin{itemize}
		\item If $\lambda \in (0,\lambda^{e,*})$, then  $0\leq \mu^{e,+}_\lambda< \mu^{e,-}_\lambda<+\infty$. 
		\item If $\lambda \in (0,\lambda^{n,*})$, then  $0\leq \mu^{n,+}_\lambda< \mu^{n,-}_\lambda<+\infty$.
	\end{itemize}
\end{cor}
\begin{proof}
	Using Lemma \ref{lemRE} and Propositions \ref{lem22}  we derive
	$$
	0\leq  \mu^{e,+}_\lambda\leq \mu_\lambda^{e,+}(u^{n,-}_\lambda)< \mu_\lambda^{e,-}(u^{n,-}_\lambda)=\mu^{e,-}_\lambda<+\infty.
	$$
	The second part  follows by the similar arguments. 
\end{proof}

	
		\begin{cor}\label{mu-compare2}
		Let $\lambda\in(0,\lambda^{e,*})$. Then  $\mu^{e,-}_\lambda<\mu^{n,-}_\lambda$.
	\end{cor}
\begin{proof}
Proposition \ref{lem22} implies the existence of $u^{n,-}_\lambda\in W^{1,2}_0$ such that $\mu^{n,-}_\lambda=\mu^{n,-}_\lambda(u^{n,-}_\lambda)$. By Lemma \ref{lemGEOM},  $\mathcal{R}^e_\lambda(t_\lambda^{e,-}(u^{n,-}_\lambda)u^{n,-}_\lambda)$ = $\mathcal{R}^n_\lambda(t_\lambda^{e,-}(u^{n,-}_\lambda)u^{n,-}_\lambda)$ and the function $t\mapsto \mathcal{R}^n_\lambda(tu^{n,-}_\lambda)$ is decreasing in $(t_\lambda^{n,-}(u^{n,-}_\lambda),t_\lambda^{e,-}(u^{n,-}_\lambda))$.
Hence 
				$$\mu_\lambda^{e,-}\leq\mathcal{R}^e_\lambda(t_\lambda^{e,-}(u^{n,-}_\lambda)u^{n,-}_\lambda)=\mathcal{R}^n_\lambda(t_\lambda^{e,-}(u^{n,-}_\lambda)u^{n,-}_\lambda)<\mathcal{R}^n_\lambda(t_\lambda^{n,-}(u^{n,-}_\lambda)u^{n,-}_\lambda)=\mu^{n,-}_\lambda.
				$$
	\end{proof}

\noindent {\bf\textit{ Conclusion of the proof of Lemma \ref{LemMain}}:} 
	
\textbf{(i)} Observe
	\begin{equation}\label{eq:LLL}
		\lambda^e(u)<\lambda^n(u),~~\forall u \in  W^{1,2}_0\setminus 0.
	\end{equation}
Indeed, since $2^{\frac{\gamma-q}{\gamma-2}}/(q\gamma^{\frac{2-q}{\gamma-2}})>1$ (see e.g., \cite{Cher}), we get
 $$
c^n_{q,\gamma}=\frac{2^{\frac{\gamma-q}{\gamma-2}}}{(q\gamma^{\frac{2-q}{\gamma-2}})}c^e_{q,\gamma}>c^e_{q,\gamma}
$$ 
which yields \eqref{eq:LLL}. Now in view of \eqref{eq:lambdaEB} we obtain \textbf{(i)}.
	
	\textbf{(ii)} By Proposition \ref{lem22}, $0<\mu^{e,-}_\lambda<+\infty$ for $\lambda \in  (0,\lambda^{e,*})$ and there exists a minimizer $ u^{e,-}_\lambda \in W^{1,2}_0\setminus 0$ of $\mu^{e,-}_\lambda(u)$ on $W^{1,2}_0\setminus 0$.
	Since  $\mu^{e,-}_\lambda( u^{e,-}_\lambda)=\mu^{e,-}_\lambda( |u^{e,-}_\lambda|)$ we my assume that $u^{e,-}_\lambda$ is nonnegative in $\Omega$. 
Applying Lemma \ref{lemRE} we obtain that $u^{e,-}_\lambda$  weakly satisfies equation \eqref{p} with $\mu=\mu^{e,-}_\lambda$ and $\Phi_{\lambda,\mu^{e,-}_\lambda}(u^{e,-}_\lambda)=0$, $\Phi_{\lambda,\mu^{e,-}_\lambda}''(u^{e,-}_\lambda)<0$. 	
	By the maximum principle  and regularity solutions for the elliptic boundary value problems \cite{trud}, it follows that $u^{e,-}_\lambda\in C^{2}(\Omega)\cap C^{1}(\overline{\Omega})$ and $u^{e,-}_\lambda>0$ in $\Omega$.

Proof of \textbf{(iii)} follows from Corollary \ref{cor: ineqStr}.

\textbf{(iv)} Notice that Lemma \ref{lemGEOM}  entails  $\mu^{n,+}_\lambda\leq \mu^{e,+}_\lambda$ for $\lambda\in(0,\lambda^{e,*})$. From here and by Corollary  \ref{mu-compare2} it follows that
							$0\leq \mu^{n,+}_\lambda\leq \mu^{e,+}_\lambda<\mu^{e,-}_\lambda<\mu^{n,-}_\lambda<+\infty$ for any $\lambda \in (0,\lambda^{e,*})$

\section{Proof of Theorem \ref{thm2}}\label{Sec:4}

%

Consider the Nehari manifold corresponding \eqref{p}
\begin{equation}\label{eq:Nehari}
	\mathcal{N}_{\lambda,\mu}=\{u \in W^{1,2}_0\setminus 0: ~\Phi_{\lambda,\mu}'(u)=0\}.
\end{equation}
Notice that since  $\mathcal{R}^n_\lambda(u)=\mu$ if, and only if, $\Phi'_{\lambda,\mu}(u)=0$, we have
\begin{equation*}\label{eq:NehariR}
		\mathcal{N}_{\lambda,\mu}=\{u \in W^{1,2}_0\setminus 0: ~\mathcal{R}^n_\lambda(u)=\mu\}.
\end{equation*}
Hence and since $\displaystyle\lim_{t\to 0}\mathcal{R}_{\lambda}^n(tu)=+\infty$ and 
	$\displaystyle\lim_{t\rightarrow\infty}\mathcal{R}_{\lambda}^n(tu)=-\infty$ (see Figure \ref{fig2}), we get that $\mathcal{N}_{\lambda,\mu} \neq \emptyset$, for any $\lambda>0$, $\mu\in \mathbb{R}$.
\begin{lem}\label{lem:coercive}
		$\Phi_{\lambda,\mu}$ is coercive on 
		$\mathcal{N}_{\lambda,\mu}$ for any $\lambda >0$, $\mu \in \mathbb{R}$. 
\end{lem}
\begin{proof}
Let $u\in \mathcal{N}_{\lambda,\mu}$. Then 
$
\|u\|_1^2+\lambda\|u\|_{L^q}^q-\mu\|u\|_{L^{\alpha}}^\alpha-\|u\|_{L^\gamma}^\gamma=0
$ and Sobolev's theorem imply
							\begin{align*}\label{Phi-bounded}
							\Phi_{\lambda,\mu}(u)=\frac{\gamma-2}{2\gamma}\|u\|_1^2+\lambda\frac{\gamma-q}{q\gamma}\|u\|_{L^q}^q-&\mu\frac{\gamma-\alpha}{\alpha\gamma}\|u\|_{L^{\alpha}}^\alpha\geq\\
							&~~~~~\frac{\gamma-2}{2\gamma}\|u\|_1^2-\mu C\|u\|_1^\alpha,
							\end{align*}
							for some $C>0$. Thus, since $\alpha<2$, $\Phi_{\lambda,\mu}(u) \to +\infty$ if  $\|u\|_1 \to +\infty$ for $u \in \mathcal{N}_{\lambda,\mu}$.
\end{proof}

Introduce, the Nehari manifold subset
\begin{equation}\label{eq:RNlm}
	\mathcal{R}\mathcal{N}_{\lambda,\mu}^+=\{u\in W^{1,2}_0\setminus 0:~ \mathcal{R}_\lambda^{n}(u)=\mu,~ (\mathcal{R}_\lambda^{n})'(u)>0\}.
\end{equation}
Observe that by Proposition \ref{Prop1},
$$
\mathcal{R}\mathcal{N}_{\lambda,\mu}^+=\{u\in W^{1,2}_0\setminus 0:~ \Phi_{\lambda,\mu}'(u)=0,~ \Phi_{\lambda,\mu}''(u)>0\}.
$$
Notice that  $\mathcal{R}\mathcal{N}_{\lambda,\mu}^+ \neq \emptyset$ for $\lambda \in (0,\lambda^{n,*})$ and $\mu \in (\mu^{n,+}_\lambda, \mu^{n,-}_\lambda)$. Indeed, since \eqref{muMinus} there is $u \in W^{1,2}_0\setminus 0$ such that $\mu^{n,+}_\lambda<\mu^{n,+}_\lambda(u) <\mu <\mu^{n,-}_\lambda<\mu^{n,-}_\lambda(u)$.  Hence the equation $\mathcal{R}^n_\lambda(tu)=\mu$ has a solution $s^1_{\lambda,\mu}(u)$ such that $\mathcal{R}^n_\lambda(s^1_{\lambda,\mu}(u)u)=\mu$ and $(\mathcal{R}^n_\lambda)'(s^1_{\lambda,\mu}(u)u)>0$. Thus $s_1(u)u \in \mathcal{R}\mathcal{N}_{\lambda,\mu}^+$.

Consider
\begin{equation} \label{GminUnst0}
	\hat{\Phi}^+_{\lambda,\mu}=\min \{\Phi_{\lambda,\mu}(u): ~u \in  \mathcal{R}\mathcal{N}_{\lambda,\mu}^+\}.
\end{equation}

\begin{prop}\label{prop:mu+}
Assume  $\lambda\in (0,\lambda^{n,*})$ and 
$\mu\in (\mu^{n,+}_\lambda, \mu^{n,-}_\lambda)$. 
Then  there exists a minimizer $u^1_{\lambda,\mu}$ of problem \eqref{GminUnst0}.

\end{prop}
\begin{proof} 
	Let $u_m$ be a minimizer sequencer of \eqref{GminUnst0}, i.e.,
	$$
	\Phi_{\lambda,\mu}(u_m) \to \hat{\Phi}^+_{\lambda,\mu},~~\Phi_{\lambda,\mu}'(u_m)=0,~ \Phi_{\lambda,\mu}''(u_m)>0,~~m=1,\ldots .
	$$
By the coerciveness of $\Phi_{\lambda,\mu}$  on $\mathcal{N}_{\lambda,\mu}$, the sequence $(u_m)$ is bounded in $W^{1,2}_0$ and therefore we may assume
$$
u_m \to \bar{u}~~\mbox{strongly in}~~L^p~~\mbox{and weakly in} ~W^{1,2}_0
$$
where $p \in(1,2^*)$.

It is easily to see that if $u_m \to \bar{u}$ strongly in $W^{1,2}_0$, then $\bar{u}$ is a minimizer  of \eqref{GminUnst0}. Suppose the converse, then $\|\bar{u}\|_1 <\liminf_{m\to \infty}\|u_m\|_1$ and 
\begin{align}
	&\Phi_{\lambda,\mu}(\bar{u})< \liminf_{m\to \infty}\Phi_{\lambda, \mu}(u_m)=\hat{\Phi}^+_{\lambda,\mu},\label{eq:PHI0} \\
	&\Phi'_{\lambda,\mu}(\bar{u})< \liminf_{m\to \infty}\Phi'_{\lambda, \mu}(u_m)=0, ~~\mathcal{R}_\lambda(\bar{u})< \liminf_{m\to \infty}\mathcal{R}_\lambda(u_m)=\mu.\label{eq:PHI00}
\end{align}
Notice that $\|u_m\|_{L^\alpha} \geq c_0$ for some $c_0>0$ which does not depend on $m=1,2,\ldots$ . Indeed, the inequality $\Phi'_{\lambda,\mu}(\bar{u})<0$ implies $\bar{u}\neq 0$. Since, $u_m\to \bar{u}$ in $L^\alpha(\Omega)$, we obtain that $\|u_m\|_{L^\alpha} \geq c_0$, for some $c_0>0$.

The inequality $\Phi'_{\lambda,\mu}(\bar{u})<0$ implies two possibilities: 1) there exist the distinct critical points $s^0_{\lambda,\mu}(\bar{u})$, $s^1_{\lambda,\mu}(\bar{u})$ such that $0<s^0_{\lambda,\mu}(\bar{u})<1<s^1_{\lambda,\mu}(\bar{u})$, or 2) $0<s^2_{\lambda,\mu}(\bar{u})<1$. 
In  case 1), we have $\Phi_{\lambda,\mu}(s^1_{\lambda,\mu}(\bar{u})\bar{u})<\Phi_{\lambda,\mu}(\bar{u})<\hat{\Phi}^+_{\lambda,\mu}$ and since $s^1_{\lambda,\mu}(\bar{u})\bar{u} \in \mathcal{R}\mathcal{N}_{\lambda,\mu}^+$, we obtain a contradiction.

Suppose 2) $s^2_{\lambda,\mu}(\bar{u})<1$. Consider the set $R_m(t):=\mathcal{R}_\lambda^n(tu_m)$, $t>0$, $m=1,2,...,$.
We claim that the sequence of the functions $R_m$ is bounded  in $C^1[\sigma, T]$ for any $\sigma, T\in (0,+\infty)$.
Indeed, for $\sigma, T\in (0,+\infty)$, $\sigma< T$, due to Sobolev's inequality $\|u\|^q_{L^q}\leq C\|u\|^q_1$ we have
	\begin{align*}
	R_m(t)\leq \frac{T^{2-\alpha}\|u_m\|_1^2+C\lambda\sigma^{q-\alpha}\|u_m\|_1^{q}}{\|u_m\|_{L^\alpha}^\alpha},~R'_m(t)\leq \frac{(2-\alpha)\sigma^{1-\alpha}\|u_m\|_1^2}{\|u_m\|_{L^\alpha}^\alpha},\label{r-m}
	\end{align*}
for  $t \in [\sigma, T]$, $m=1,2,\ldots$, where $0<C<+\infty$ does not depend on $m=1,2,\ldots$. Hence, due to  $\|u_m\|_{L^\alpha} \geq c_0$ for $m=1,2,...$ and by the boundedness of $\|u_m\|_1$,  we get that $(R_m)$ is bounded  in $C^1[\sigma, T]$.

Thus by the Arzela-Ascoli compactness criterion
we can assume that for any $\sigma, T\in (0,+\infty)$, there holds
\begin{equation}\label{ravConv}
	R_m(t) \to \bar{R}(t) ~~\mbox{in}~~C[\sigma, T]~~\mbox{as}~~m\to \infty
\end{equation}
for some limit function $\bar{R} \in C(0, +\infty)$.
Evidently the sequences $(s^0_{\lambda,\mu}(u_m))$, $(s^1_{\lambda,\mu}(u_m))$ are bounded. The sequence $s^2_{\lambda,\mu}(u_m)$ is also bounded, since $(u_m)$ is separated from zero and bounded in $W^{1,2}_0$. Indeed, in the converse case, we get a contradiction
\begin{align*}
	\mu=\mathcal{R}^n_\lambda(&s^2_{\lambda,\mu}(u_m)u_m):=
	\frac{1}{\|u_m\|_{L^{\alpha}}^\alpha}\times ((s^2_{\lambda,\mu}(u_m))^{2-\alpha}\| u_m\|^2_1+\\&
	(s^2_{\lambda,\mu}(u_m))^{q-\alpha}\lambda\|u_m\|^q_{L^q}-(s^2_{\lambda,\mu}(u_m))^{\gamma-\alpha}\|u_m\|^\gamma_{L^\gamma})\to -\infty.
\end{align*}
Hence we may assume that
$$
s^0_{\lambda,\mu}(u_m) \to \bar{s}^0,~~s^1_{\lambda,\mu}(u_m) \equiv 1=\bar{s}^1,~~s^2_{\lambda,\mu}(u_m) \to \bar{s}^2~~\mbox{as}~~m \to +\infty.
$$
It easy to see that $\bar{R}(\bar{s}^0)=\bar{R}(\bar{s}^1)=\bar{R}(\bar{s}^2)=\mu$, $0<\bar{s}^0\leq\bar{s}^1\leq \bar{s}^2<+\infty$ and 
\begin{equation}\label{cases0}
	\left\{
	\begin{aligned}
		&\bar{R}(s)\geq\mu~~\mbox{if}~~ s \in (0,\bar{s}^0),\\
		&\bar{R}(s)\leq \mu~~\mbox{if}~~ s \in [\bar{s}^0, \bar{s}^1],\\
		&\bar{R}(s) \geq\mu~~\mbox{if}~~ s \in (\bar{s}^1,\bar{s}^2).
	\end{aligned}
	\right.
\end{equation}
From the above it follows that $\mathcal{R}^n_\lambda(s\bar{u})<\bar{R}(s)$ for $s>0$. Furthermore, \eqref{mumu}  implies $\mu<\mu^{n,-}_\lambda=\inf_{W^{1,2}_0\setminus 0}\mu^{n,-}_\lambda(u)\leq \mathcal{R}^n_\lambda(t_\lambda^{n,-}(\bar{u})\bar{u})$. Hence, in view of \eqref{cases0}, there is only the following two  possibilities:
a)  $t_\lambda^{n,-}(\bar{u}) \in (\bar{s}^1,\bar{s}^2)$; ~ b) $t_\lambda^{n,-}(\bar{u}) \in (0,\bar{s}^0)$.

In case a), we get a contradiction $1=s^1_{\lambda,\mu}(\bar{u})<s^2_{\lambda,\mu}(\bar{u})<1$ . 

Suppose b). Observe
$$
0<(\mathcal{R}^n_\lambda)'(s^1_{\lambda,\mu}(\bar{u})\bar{u})< \liminf_{j\to \infty}(\mathcal{R}^n_\lambda)'(s^1_{\lambda,\mu}(\bar{u}))u_m).
$$
Hence,  for sufficiently large $m$, $(\mathcal{R}^n_\lambda)'(s^1_{\lambda,\mu}(\bar{u}))u_m)>0$  and therefore $s^1_{\lambda,\mu}(\bar{u})>t_\lambda^{n,+}(u_m)>s^0_{\lambda,\mu}(u_m)$. However,  $t_\lambda^{n,-}(\bar{u}) \in (0,\bar{s}^0)$ entails $s^1_{\lambda,\mu}(\bar{u})<\bar{s}^0$ and thus $s^1_{\lambda,\mu}(\bar{u})<s^0_{\lambda,\mu}(u_m)$ for sufficiently large $m$  which  contradicts to the inequality $(\mathcal{R}^n_\lambda)'(s^1_{\lambda,\mu}(\bar{u}))u_m)>0$. 
This concludes the proof.
	\end{proof}

\medskip

\textit{\textbf{Conclusion of the proof of  Theorem \ref{thm2}}.}

 Let $\lambda\in (0,\lambda^{e,*})$. Since $\lambda^{e,*}<\lambda^{n,*}$, by Proposition \ref{prop:mu+}  there exists a minimizer $u^1_{\lambda,\mu}$ of problem \eqref{GminUnst0} for any $\mu\in (\mu^{n,+}_\lambda, \mu^{n,-}_\lambda)$. Since $\Phi_{\lambda, \mu}(u^1_{\lambda, \mu})=\Phi_{\lambda, \mu}(|u^1_{\lambda, \mu}|)$, $\Phi_{\lambda, \mu}'(u^1_{\lambda, \mu})=\Phi_{\lambda, \mu}'(|u^1_{\lambda, \mu}|)$, $\Phi_{\lambda, \mu}''(u^1_{\lambda, \mu})=\Phi_{\lambda, \mu}''(|u^1_{\lambda, \mu}|)$ we may assume that $u^1_{\lambda,\mu}$ is a non-negative function on $\Omega$. 

Let us prove that 
\begin{equation}\label{cases}
\hat{\Phi}^+_{\lambda,\mu}\equiv \Phi_{\lambda, \mu}(u^1_{\lambda, \mu})<0~~\mbox{if}~~\mu  \in (\mu^{e,+}_\lambda, \mu^{n,-}_\lambda),~~\lambda \in (0,\lambda^{e,*}). 
\end{equation}	
Since $\lambda \in (0,\lambda^{e,*})$,  $\mu^{e,+}_\lambda< \mu^{n,-}_\lambda$ and by the definitions of $\mu^{e,+}_\lambda, \mu^{n,-}_\lambda$  there is $\hat{u} \in W^{1,2}_0\setminus 0$ such that $\mu^{n,-}_\lambda(\hat{u})>\mu>\mu^{e,+}_\lambda(\hat{u})$.  Then  $s_1(\hat{u}) \in (t^{e,+}_\lambda(\hat{u}), t^{n,-}_\lambda(\hat{u})) \subset  (t^{e,+}_\lambda(\hat{u}), t^{e,-}_\lambda(\hat{u}))$. Since $\mathcal{R}_\lambda^n(t\hat{u})>\mathcal{R}_\lambda^e(t\hat{u})$ for $t \in (t^{e,+}_\lambda(\hat{u}), t^{e,-}_\lambda(\hat{u}))$, we obtain that 
$\mu=\mathcal{R}_\lambda^n(s_1(\hat{u})\hat{u})>\mathcal{R}_\lambda^e(s_1(\hat{u})\hat{u})$, which implies $\Phi_{\lambda, \mu}(s_1(\hat{u})\hat{u}) <0$ and therefore
$$\Phi_{\lambda, \mu}(u^1_{\lambda, \mu}) =\hat{\Phi}^+_{\lambda,\mu} \leq \Phi_{\lambda, \mu}(s_1(\hat{u})\hat{u}) <0,
$$
since $s_1(\hat{u})\hat{u} \in \mathcal{R}\mathcal{N}_{\lambda,\mu}^+$. 
Thus we get \eqref{cases}.

To show that $u^1_{\lambda,\mu}$ is weakly satisfies  to \eqref{p}, it suffices to have that $\Phi_{\lambda, \mu}''(u^1_{\lambda, \mu})>0$. Indeed, in this case, by the Lagrange multiplier rules there exist $\nu_0, \nu_1$ such that $|\nu_0|+|\nu_1|\neq 0$ and
$$
\nu_0 D_u\Phi_{\lambda, \mu}(u^1_{\lambda, \mu})+\nu_1 D_u\Phi'_{\lambda, \mu}(u^1_{\lambda, \mu})=0.
$$
Testing this equality by $u^1_{\lambda, \mu}$ we obtain $\nu_1\Phi''_{\lambda, \mu}(u^1_{\lambda, \mu})=0$ which implies that $\nu_1=0$ and 
consequently we get the desired conclusion. 

To prove  $\Phi_{\lambda, \mu}''(u^1_{\lambda, \mu})>0$ for  $\mu \in (\mu^{e,+}_\lambda,\mu^{n,-}_\lambda)$, it is sufficient to show that the strong inequalities  $s^0_{\lambda,\mu}(u^1_{\lambda, \mu})<s^1_{\lambda,\mu}(u^1_{\lambda, \mu})<s^2_{\lambda,\mu}(u^1_{\lambda, \mu})$ hold. 
Notice that by the construction we have $u^1_{\lambda, \mu}\in \mathcal{R}\mathcal{N}_{\lambda,\mu}^+$ and thus $s^1_{\lambda,\mu}(u^1_{\lambda, \mu})$ is well defined so that $s^0_{\lambda,\mu}(u^1_{\lambda, \mu})\leq s^1_{\lambda,\mu}(u^1_{\lambda, \mu})$. 
Since $\mu <\mu^{n,-}_\lambda\leq \mu^{n,-}_\lambda(u^1_{\lambda, \mu})$, Proposition \ref{p28} implies  $s^1_{\lambda,\mu}(u^1_{\lambda, \mu})<s^2_{\lambda,\mu}(u^1_{\lambda, \mu})$.  Applying Lemma \ref{lemGEOM} (see Figure \ref{fig4}) we derive $\mathcal{R}^e_{\lambda}(s^0_{\lambda,\mu}(u^1_{\lambda, \mu})u)>\mathcal{R}^n_{\lambda}(s^0_{\lambda,\mu}(u^1_{\lambda, \mu})u)=\mu$. Consequently, $\Phi_{\lambda, \mu}(s^0_{\lambda,\mu}(u^1_{\lambda, \mu})u^1_{\lambda, \mu})>0$ for $\mu \in (\mu^{e,+}_\lambda,\mu^{n,-}_\lambda)$. On the other hand, by \eqref{cases},  $\Phi_{\lambda, \mu}(s^1_{\lambda,\mu}(u^1_{\lambda, \mu})u^1_{\lambda, \mu})<0$ and therefore $s^0_{\lambda,\mu}(u^1_{\lambda, \mu})<s^1_{\lambda,\mu}(u^1_{\lambda, \mu})$ for any $\mu \in (\mu^{e,+}_\lambda,\mu^{n,-}_\lambda)$. Thus, indeed, $\Phi_{\lambda, \mu}''(u^1_{\lambda, \mu})>0$ and consequently  $u^1_{\lambda, \mu}$ is a weak non-negative solution of \eqref{p}, for any $\mu \in (\mu^{e,+}_\lambda,\mu^{n,-}_\lambda)$.

By the maximum principle  and the regularity solutions of elliptic boundary value problems \cite{trud} it follows that $u^1_{\lambda, \mu}\in C^{2}(\Omega)\cap C^{1}(\overline{\Omega})$ and $u^1_{\lambda, \mu}(x)>0$ in $\Omega$.  This concludes  the proof of  Theorem \ref{thm2}.

\section{Proof of Theorem \ref{thm3}}\label{Sec:5}

Consider the following subset of the Nehari manifold $\mathcal{N}_{\lambda,\mu}$
\begin{equation}\label{eq:RNlm-}
	\mathcal{R}\mathcal{N}_{\lambda,\mu}^{2}=\{u\in W^{1,2}_0\setminus 0:~ \mathcal{R}_{\lambda}^n(u)=\mu,~ \mathcal{R}_\lambda'(u)<0,~~\mathcal{R}_\lambda''(u)<0\}.
\end{equation}
Observe that by Proposition \ref{Prop2},
$$
\mathcal{R}\mathcal{N}_{\lambda,\mu}^{2}=\{u\in W^{1,2}_0\setminus 0:~ \Phi_{\lambda,\mu}'(u)=0,~ \Phi_{\lambda,\mu}''(u)<0~~\mathcal{R}_\lambda''(u)<0\}.
$$
The set $ \mathcal{R}\mathcal{N}_{\lambda,\mu}^{2}$ is not empty if $\lambda\in (0,\lambda^{n,*})$ and $\mu\leq \mu^{n,-}_\lambda$. Indeed, let $u\in W^{1,2}_0\setminus 0$. Then, since $\lambda\in (0,\lambda^{n,*})$, Proposition \ref{propRN2} implies that the function $\mathcal{R}^n_\lambda(tu)$ has precisely two distinct critical points $t_\lambda^{n,+}(u),t_\lambda^{n,-}(u)$. Notice that $\displaystyle\lim_{t\rightarrow\infty}\mathcal{R}_{\lambda}^n(tu)=-\infty$ and $t\mapsto \mathcal{R}_{\lambda}^n(tu)$ is decreasing function in the interval $(t_\lambda^{n,-}(u),\infty)$ for any $u\in W^{1,2}_0\setminus 0$. Therefore and since $\mu <\mu^{n,-}_\lambda\leq \mu_\lambda^{n,-}(u)$, there exists a unique $s^2_{\lambda, \mu}(u)>t_\lambda^{n,-}(u)$ such that $\mu=\mathcal{R}_{\lambda}^n(s^2_{\lambda, \mu}(u)u)$ and $(\mathcal{R}_{\lambda}^n)'(s^2_{\lambda, \mu}(u)u)<0$.
Proposition \ref{prop:2deriv} implies that there is $r^n_\lambda(u)$ such that $t_\lambda^{n,+}(u)<r^n_\lambda(u)<t_\lambda^{n,-}(u)$ and $(\mathcal{R}^n_\lambda)''(tu)<0$ for $t \in (r_\lambda, +\infty)$. Since $r^n_\lambda(u)<t_\lambda^{n,-}(u)<s^2_{\lambda, \mu}(u)$, we have $(\mathcal{R}^n_\lambda)''(s^2_{\lambda, \mu}(u)u)<0$.
Hence $s^2_{\lambda, \mu}(u)u \in \mathcal{R}\mathcal{N}_{\lambda,\mu}^{2}$.

Consider
\begin{equation} \label{GminUnst1}
	\hat{\Phi}^{2}_{\lambda,\mu}=\min \{\Phi_{\lambda,\mu}(u): ~u \in  \mathcal{R}\mathcal{N}_{\lambda,\mu}^{2}\}.
\end{equation}

%
%

\begin{prop}\label{prop:mu-}
Assume  $\lambda\in (0,\lambda^{n,*})$ and $\mu\leq \mu^{n,-}_\lambda$. 
Then there exists a minimizer $u^2_{\lambda,\mu}$ of problem \eqref{GminUnst1}.
\end{prop}

\begin{proof} There is $C>0$ which does not depend on $u \in\mathcal{R}\mathcal{N}_{\lambda,\mu}^{2}$ such that $\|s^2_{\lambda, \mu}(u)u\|_{L^\gamma} >C$ for all $u\in  W^{1,2}_0\setminus 0$. Indeed, by Lemma \ref{lemGEOM} and due to embedding $W^{1,2}_0\hookrightarrow L^\gamma (\Omega)$,  we have
$$
s^2_{\lambda, \mu}(u)>t_\lambda^{n,-}(u)>t^{n}(u)=C^n(\alpha,q,\gamma)\left(\frac{\|u\|_1^2}{\|u\|_{L^\gamma}^{\gamma}}\right)^{1/(\gamma-2)}\geq \frac{C}{\|u\|_{L^\gamma}}.
$$
Thus,	 
\begin{equation}
\label{n73}
s^2_{\lambda, \mu}(u)\|u\|_{L^\gamma}\geq C>0.
\end{equation}
Let $(u_m)$ be a minimizer sequence of \eqref{GminUnst1}, that is
					$$
					\Phi_{\lambda,\mu}(u_m) \to \hat{\Phi}^{2}_{\lambda,\mu},~\mbox{as}~m\to \infty,~~\Phi_{\lambda,\mu}'(u_m)=0,~ \Phi_{\lambda,\mu}''(u_m)<0,~~\mathcal{R}_\lambda''(u_m)<0.
					$$	
Then
\begin{equation}
\label{n74}
s^2_{\lambda, \mu}(u_m)=1,~\mbox{for}~m=1,2,\cdots.
\end{equation}	
By Lemma \ref{lem:coercive}, 
$\Phi_{\lambda,\mu}$ is coercive on $\mathcal{N}_{\lambda,\mu}$ and thus $(u_m)$ is bounded in $W^{1,2}_0$. Hence by the Sobolev embedding  and  Banach-Alaoglu  theorems, up to a subsequence, there holds 
					$$
					u_m \to \hat{u}~~\mbox{strongly in}~~L^p~~\mbox{and weakly in} ~W^{1,2}_0
					$$
					where $1< p<2^*$.
In particular, we obtain from \eqref{n73} and \eqref{n74} that
 $\hat{u}\neq 0$.

 Let us show  that, up to a subsequence, there holds a strong convergence  $u_m \to \hat{u}$ in $W^{1,2}_0$. Otherwise, we would have 
					\begin{align}  \label{eq:RL3}
					&\mathcal{R}^n_\lambda(\hat{u})< \liminf_{m\to \infty}\mathcal{R}^n_\lambda(u_m)=\mu,\\ \label{eq:RL31}
					&(\mathcal{R}^n_\lambda)'(\hat{u})< \liminf_{m\to \infty}(\mathcal{R}^n_\lambda)'(u_m)= 0,\\ \label{eq:RL}
					&(\mathcal{R}^n_\lambda)''(\hat{u})< \liminf_{m\to \infty}(\mathcal{R}^n_\lambda)''(u_m) \leq 0,				
\end{align}
and
\begin{equation}\label{ineqPHI}
	0=\Phi'_{\lambda,\mu}(s^2_{\lambda, \mu}(\hat{u})\hat{u})< \liminf_{m\to \infty}\Phi'_{\lambda,\mu}(s^2_{\lambda, \mu}(\hat{u})u_m).
\end{equation}
In particular, the last inequality implies that 
$\Phi'_{\lambda,\mu}(s^2_{\lambda, \mu}(\hat{u})u_m)>0$ holds for sufficiently large $m$.

Since $s^2_{\lambda, \mu}(\hat{u})>t_\lambda^{n,-}(\hat{u})$ and $\mu=\mathcal{R}_{\lambda}^n(s^2_{\lambda, \mu}(\hat{u})\hat{u})$ and $(\mathcal{R}_{\lambda}^n)'(s^2_{\lambda, \mu}(\hat{u})\hat{u})<0$, $(\mathcal{R}_{\lambda}^n)''(s^2_{\lambda, \mu}(\hat{u})\hat{u})<0$, we get from  
 \eqref{eq:RL3}-\eqref{eq:RL} and Proposition \ref{prop:2deriv} the strong inequality  $s^2_{\lambda, \mu}(\hat{u})<1$.  By setting $R_m(t):=\mathcal{R}_\lambda^n(tu_m)$, $t>0$, $m=1,2,...,$ and proceeding as done in the proof of Theorem \ref{thm2}, we have that  $R_m(t) \to \bar{R}(t)$ in $C[\sigma, T]$,  for some limit function $\bar{R} \in C(0, +\infty)$ and  any $\sigma, T\in (0,+\infty)$. Furthermore, there  holds 
\begin{equation}\label{convergM}
	s^0_{\lambda,\mu}(u_m) \to \bar{s}^0,~~s^1_{\lambda,\mu}(u_m) \to \bar{s}^1,~~s^2_{\lambda,\mu}(u_m) =\bar{s}^2=1~~\mbox{as}~~m \to +\infty.
\end{equation}
so that  $\bar{R}(\bar{s}^0)=\bar{R}(\bar{s}^1)=\bar{R}(\bar{s}^2)=\mu$ and $\bar{R}(t)$ satisfies to \eqref{cases0}. 

As done in the proof of Theorem \ref{thm2},we have just  two  possibilities:
 $a)$ $t_\lambda^{n,-}(\bar{u}) \in (0,\bar{s}^0)$  and $b)$ $t_\lambda^{n,-}(\bar{u}) \in (\bar{s}^1,\bar{s}^2=1)$. Suppose  $a)$. Then $s^2_{\lambda, \mu}(\hat{u}) \in (0,  \bar{s}^0)$ and consequently by \eqref{convergM} we get  $0<s^2_{\lambda, \mu}(\hat{u})<s^0_{\lambda, \mu}(u_m)$ for sufficiently large $m$.  However, this implies
that $\Phi'_{\lambda, \mu}(s^2_{\lambda, \mu}(\hat{u})u_m)<0$ for sufficiently large $m$, which contradicts to \eqref{ineqPHI}.

Assume $ b)$. Then, due to \eqref{convergM}, we have $s^1_{\lambda, \mu}(u_m)< t_\lambda^{n,-}(\hat{u}) 
<s^2_{\lambda, \mu}(\hat{u})<1=s^2_{\lambda, \mu}(u_m)$ for sufficiently large $m$. Besides this, we know that  $\Phi'_{\lambda, \mu}(tu_m)>0$ for any $s^1_{\lambda, \mu}(u_m)<t<s^2_{\lambda, \mu}(u_m)$. Thus, we have by contradiction assumption, that
$$
\Phi_{\lambda,\mu}(s^2_{\lambda, \mu}(\hat{u})\hat{u})<\liminf_{m\to \infty}\Phi_{\lambda,\mu}(s^2_{\lambda, \mu}(\hat{u})u_m)\leq \liminf_{m\to \infty}\Phi_{\lambda,\mu}(u_m)=\hat{\Phi}^{2}_{\lambda,\mu},
$$
which implies a contradiction because  $(s^2_{\lambda, \mu}(\hat{u})\hat{u} \in \mathcal{R}\mathcal{N}_{\lambda,\mu}^2$.

Thus $u_m \to \hat{u}$ strongly in $W^{1,2}_0$. To complete the proof of Proposition, we just remain to show that $u^2_{\lambda,\mu} \in \mathcal{R}\mathcal{N}_{\lambda,\mu}^2$.
Since 
$0<\lambda <\lambda^{n,*}$, Proposition \ref{propRN2} yields that the strong inequality $0<t_\lambda^{n,+}(u^2_{\lambda,\mu})<t_\lambda^{n,-}(u^2_{\lambda,\mu})$ holds. So, there exists  $r_\lambda(u^2_{\lambda,\mu}) \in (t_\lambda^{n,+}(u^2_{\lambda,\mu}), t_\lambda^{n,-}(u^2_{\lambda,\mu}))$ such that $(\mathcal{R}_\lambda^n)'(r_\lambda(u^2_{\lambda,\mu}) u^2_{\lambda,\mu})>0$ and $(\mathcal{R}_\lambda^n)^{\prime\prime}(r_\lambda(u^2_{\lambda,\mu}) u^2_{\lambda,\mu})=0$. Since $(\mathcal{R}_\lambda^n)'(u^2_{\lambda,\mu})\leq 0$ (due to strong convergence), it follows from Proposition \ref{prop:2deriv}, that $r_\lambda(u^2_{\lambda,\mu})<1$, that is, $(\mathcal{R}_\lambda^n)^{\prime\prime}( u^2_{\lambda,\mu})<0$. In particular, this inequality implies that $t_\lambda^{n,-}(u^2_{\lambda,\mu}) <1$, that is, $(\mathcal{R}_\lambda^n)^{\prime} (u^2_{\lambda,\mu})<0$. This implies that $u^2_{\lambda,\mu} \in \mathcal{R}\mathcal{N}_{\lambda,\mu}^2$.

Finally, we note that  the fact that $u_m \to \hat{u}$ strongly in $W^{1,2}_0$ and $u^2_{\lambda,\mu} \in \mathcal{R}\mathcal{N}_{\lambda,\mu}^2$ yield that $\hat{u}$ is a minimizer of \eqref{GminUnst1}.   This ends the proof.
\end{proof}
\medskip

\medskip

\noindent\textit{\textbf{Conclusion of the proof of  Theorem $\ref{thm3}$}.} Assume  $0<\lambda<\lambda^{n,*}$ and $-\infty<\mu\leq \mu^{n,-}_\lambda$. Then,  by Proposition \ref{prop:mu-}, there exists a minimizer $u^2_{\lambda,\mu}\in W^{1,2}_0\setminus 0$ of problem \eqref{GminUnst1}. Since $\Phi_{\lambda,\mu}''(u^2_{\lambda,\mu})<0$, we obtain by similar arguments as done in the proof  of Theorem  \ref{thm2} that  $u^2_{\lambda,\mu}$ is a weak positive  $C^{2}(\Omega)\cap C^{1}(\overline{\Omega})$-solution of problem \eqref{p}.  

To finish the proof, let us show \eqref{Phi-}. First, assume $ \mu^{e,-}_\lambda<\mu<\mu^{n,-}_\lambda$. Then there exists  $u \in W^{1,2}_0 \setminus 0$ such that $\mu^{e,-}_\lambda<\mu^{e,-}_\lambda(u)<\mu<\mu^{n,-}_\lambda<\mu^{n,-}_\lambda(u)$, that is,  $t^{n,-}_{\lambda,\mu}(u)<s^2_{\lambda,\mu}(u)<t^{e,-}_{\lambda,\mu}(u)$ with $s^2_{\lambda,\mu}(u)u\in \mathcal{R}\mathcal{N}_{\lambda,\mu}^{2}$. This and Lemma \ref{lemGEOM} imply 
$ \mathcal{R}^e_{\lambda}(s^2_{\lambda,\mu}(u)u)< \mathcal{R}^n_{\lambda}(s^2_{\lambda,\mu}(u)u)=\mu$ and consequently $\Phi_{\lambda,\mu}(s^2_{\lambda,\mu}(u)u)<0$. Since $s^2_{\lambda,\mu}(u)u\in \mathcal{R}\mathcal{N}_{\lambda,\mu}^{2}$, we obtain from \eqref{GminUnst1} that 
$\Phi_{\lambda,\mu}(u^2_{\lambda, \mu})\leq \Phi_{\lambda,\mu}(s^2_{\lambda,\mu}(u)u)<0$.

If $\mu=\mu^{e,-}_\lambda$, then $\mathcal{R}^e_{\lambda}(u^2_{\lambda, \mu^{e,-}}) = \mathcal{R}^n_{\lambda}(u^2_{\lambda, \mu^{e,-}})=\mu^{e,-}$, and therefore we have $\Phi_{\lambda,\mu^{e,-}_\lambda}(u^2_{\lambda, \mu^{e,-}_\lambda})=0$. Finally, let  $ -\infty<\mu< \mu^{e,-}_\lambda$. Since  $\mathcal{R}_\lambda^n( u^2_{\lambda,\mu})=\mu$, we have that $1> t_\lambda^{e,-}(u^2_{\lambda, \mu})$ and thus Lemma \ref{lemGEOM} implies that
$\mathcal{R}^e_{\lambda}(u^2_{\lambda, \mu}) > \mathcal{R}^n_{\lambda}(u^2_{\lambda, \mu})=\mu$
or equivalently $
\Phi_{\lambda,\mu}(u^2_{\lambda, \mu})>0,~~~\mu \in (-\infty, \mu^{e,-}_\lambda).
$
This ends the proof of Theorem.

\section{Appendix}\label{Sec:appendix}

In this Appendix,  we present some results  in general sense.
Assume $W$ is a real Banach space,  $T(u)$ and $G(u)$ are   Fr\'{e}chet-differentiable functionals with derivatives $D_u T(u)$, $D_uG(u)$. 
 We deal with the equations in the following abstract form
\begin{equation}
\label{I}
D_u\Phi_\nu(u)\equiv D_u T(u)-\nu D_uG(u)=0,~~~u \in W, 
\end{equation}
where $\Phi_\nu:=T(u)-\nu G(u)$ is a corresponding variational (energy) functional and $\nu \in \mathbb{R} $ is a parameter. 
We assume that  $D_uG(u)(u) > 0$,  $\forall u\in W \setminus 0$. The case $D_uG(u)(u) < 0$, $\forall u\in W \setminus 0$ is considered analogously.  Let us also suppose that  $D_uT(tu)(tu)$, $D_uG(tu)(tu): \mathbb{R}^+ \times (W\setminus 0) \to \mathbb{R}$ are maps of $C^k$-class with $k\geq 1$.

The Nehari manifold associated with  (\ref{I}) is defined as follows 
$$
\mathcal{N}_\nu:=\{u \in W\setminus 0:~\Phi'_\nu(u):= D_u\Phi_\nu(u)(u)=0\}.
$$
We already know that $\Phi''_\nu(u)\neq 0$ for  some $u \in N_\nu$ implies that  $N_\nu$ is a $C^1$-manifold of codimension $1$, in some neighbourhood of $u$,  so that $W=T_{u}(\mathcal{N}_\nu) \oplus \mathbb{R}u$ (see \cite{ilyaReil}). 

Consider the Rayleigh quotients  
$$
\mathcal{R}^n(u):=\frac{T'(u)}{G'(u)}~~\mbox{and}~~\mathcal{R}^e(u):=\frac{T(u)}{G(u)}.
$$
Here, as before,  we are denoting by $T'(u):=d T(tu)/dt|_{t=1}=D_u T(u)(u)$ and $G'(u):=d G(tu)/dt|_{t=1}=D_uG(u)(u)$. In what follows, we suppose that $G(u),G'(u) \neq 0$ for every $u\in W\setminus 0$.

After these, we are able to define  the fibering functions 
$$\mathcal{R}^n(tu):=T'(tu)/G'(tu)~~\mbox{and}~~\mathcal{R}^e(tu):=T(tu)/G(tu),~t\geq 0,$$  defined 
for any $u \in W\setminus 0$ given. When $(\mathcal{R}^n)'(t_0u)=0$, for some $t_0>0$, we say that $t_0u$ is a \textit{fibering critical point} of $\mathcal{R}^n(tu)$; otherwise, $t_0u$ is called a \textit{fibering regular point}. 

Let $U$ be an open subset of $W\setminus 0$.  A point $\nu \in \mathbb{R}$ is said to be a \textit{fibering regular value} of $\mathcal{R}^n(u)$ on $U$ if all points $tu$, with $u \in U$ and $t \in \mathbb{R}$, in the pre-image $(\mathcal{R}^n)^{-1}(\nu)$ are fibering regular points. 
We call $t_0 \in \mathbb{R}^+$ the \textit{extreme point} of $\mathcal{R}^n(tu)$ if the function $\mathcal{R}^n(tu)$ achieves at $t_0$ its local maximum or minimum on $\mathbb{R}^+$.

As a consequence of the fact that $tu$ belongs to $\mathcal{N}_\nu$ if, and only if, it lies on the level set $\mathcal{R}^n(tu)=\nu$, we obtain
\begin{equation}\label{Req1}
(\mathcal{R}^n)'(tu) =
\frac{1}{G'(tu)} \Phi''_{\nu}(tu), ~\forall ~u\in \mathcal{N}_\nu,
\end{equation}
which implies the next result
\begin{prop} \label{Prop1}
Assume  $\mathcal{R}^n(tu)=\nu$. Then:
\begin{enumerate}
\item[$(i)$] $(\mathcal{R}^n)'(tu)=0$   if, and only if, $\Phi''_\nu(tu)=0$,
\item[$(ii)$] $(\mathcal{R}^n)'(tu)>0$  if, and only if, $\Phi''_\nu(tu)>0$,
\item[$(iii)$] $(\mathcal{R}^n)'(tu)<0$   if, and only if, $\Phi''_\nu(tu)<0$.
\end{enumerate}
\end{prop}

After these, we have the following result (see e.g. \cite{ilyaReil}).
\begin{lem} \label{LP}
	 Let $U$ be an open subset of $W\setminus 0$. Assume that  $\nu \in \mathbb{R}$ is a fibered regular value of $\mathcal{R}^n(tu)$ on $U$, then  $N_\nu\cap U$ is a $C^1$-manifold of codimension $1$ so that $W=T_{u}(\mathcal{N}_\nu) \oplus \mathbb{R}u$ for every $u \in N_\nu\cap U$.
\end{lem}

Finally, by letting $\mathcal{R}^e(tu)=\nu$, we obtain
\begin{equation}
	(\mathcal{R}^e(tu))' =\frac{1}{G(tu)}(T'(tu)-\nu G'(tu))=\frac{1}{G(tu)}\Phi'_{\nu}(tu),
\end{equation}
and 
	\begin{align*}
		(\mathcal{R}^e(tu))'' &=\frac{1}{G(tu)}(T''(tu)-\nu G''(tu))-2\frac{G'(tu)}{G^2(tu)}(T'(tu)-\nu G'(tu))\\
		&=\frac{1}{G'(tu)} \Phi''_{\nu}(tu)-2\frac{G'(tu)}{G^2(tu)}\Phi'_{\nu}(tu),
	\end{align*}
which implies 
\begin{prop} \label{Prop2}  If  $\mathcal{R}^e(tu)=\nu$ (or the same $\Phi_\nu(tu)=0$), then:
\begin{enumerate}
\item[$(i)$] $D_u\mathcal{R}^e(u)=0$ if, and only if, $D_u\Phi_\nu(u)=0$,
\item[$(ii)$] $(\mathcal{R}^e(tu))'=0$ if, and only if, $\Phi'_\nu(tu)=0$,
\item[$(iii)$] $(\mathcal{R}^e(tu))''=0$ and $(\mathcal{R}^e(tu))'=0$ if, and only if, $\Phi''_\nu(tu)=0$ and  $\Phi'_\nu(tu)=0$,
\item[$(iv)$] $(\mathcal{R}^e(tu))''>0$ and $(\mathcal{R}^e(tu))'=0$  if, and only if, $\Phi''_\nu(tu)>0$ and $\Phi'_\nu(tu)=0$,
\item[$(v)$] $(\mathcal{R}^e(tu))''<0$ and $(\mathcal{R}^e(tu))'=0$   if, and only if, $\Phi''_\nu(tu)<0$ and  $\Phi'_\nu(tu)=0$.
\end{enumerate}
\end{prop}

\end{document}